\documentclass[11pt,leqno]{article}

\usepackage[utf8]{inputenc}
\usepackage[T1]{fontenc}
\usepackage{lmodern}
\usepackage[margin=2.5cm]{geometry}  
\usepackage{graphicx}
\usepackage{tikz}
\usetikzlibrary{arrows.meta, positioning, calc}
\usepackage{xcolor}
\usepackage{amsmath,amssymb,mathtools}
\usepackage{amsthm}
\usepackage{enumitem}
\usepackage{hyperref}
\hypersetup{colorlinks=true, citecolor=teal, linkcolor=black, urlcolor=blue}
\usepackage[capitalise]{cleveref}
\usepackage[backend=biber, style=alphabetic]{biblatex}
\allowdisplaybreaks
\numberwithin{equation}{section}

\newcommand{\R}{\mathbb{R}}
\newcommand{\C}{\mathbb{C}}
\newcommand{\N}{\mathbb{N}}
\newcommand{\Z}{\mathbb{Z}}

\providecommand{\Sp}{\mathcal{S}}
\providecommand{\eps}{\varepsilon}
\providecommand{\cA}{\mathcal{A}}
\providecommand{\cB}{\mathcal{B}}
\providecommand{\cH}{\mathcal{H}}
\providecommand{\cM}{\mathcal{M}}

\providecommand{\Dbar}{\overline{\mathbb{D}}}
\providecommand{\sgn}{\operatorname{sign}}

\theoremstyle{plain}
\newtheorem{theorem}{Theorem}[section]
\newtheorem{lemma}{Lemma}[section]
\newtheorem{proposition}{Proposition}[section]
\newtheorem{corollary}{Corollary}[section]
\theoremstyle{definition}
\newtheorem{definition}{Definition}[section]
\theoremstyle{remark}
\newtheorem{remark}{Remark}[section]

\newenvironment{keyword}{\par\smallskip\noindent\textbf{Keywords: }}{\par\smallskip}

\begin{document}

\title{
Effect of a fourth-order diffusion term on the robustness of the transfer function stability
for hyperbolic systems}

\author{
  Eduardo Cerpa\footnote{Instituto de Ingenier\'ia Matem\'atica y Computacional, Facultad de Matem\'aticas, Pontificia Universidad Católica de Chile, Avda.
Vicu\~na Mackenna 4860, Macul, Santiago, Chile. E-mail: { \tt eduardo.cerpa@uc.cl}}
  \and 
  Felipe Labra\footnote{Facultad de Matem\'aticas, Pontificia Universidad Católica de Chile, Avda.
Vicu\~na Mackenna 4860, Macul, Santiago, Chile. E-mail: { \tt flabrai@estudiante.uc.cl}}
}
\date{November 4, 2026}
\maketitle


\begin{abstract}
We study a $2\times2$ hyperbolic system with anti-located control and observation. It is well known that the system
can be stabilized with a delayed output-feedback control and that, unfortunately, this design is arbitrarily sensitive to small errors in
the transport velocity. It was recently proved that adding a second-order diffusion term to the system restores
robustness. We show that the same phenomenon occurs with a fourth-order diffusion term: the maximal spectral abscissa of
the closed-loop system converges, as the diffusion parameter goes to zero, to the explicit value $-\ln(2)$, and the
resulting stability margin is robust to small perturbations of the transport velocity. This additional term leads to a
quartic characteristic equation whose roots have no tractable closed-form expression. We develop an asymptotic
classification of its four roots and use it to obtain a first-order expansion of the transfer function with an explicit
error bound. Combining this expansion with a topological degree argument, we show that $-\ln(2)$ is in fact the exact
limiting value of the spectral abscissa, extending to the fourth-order case the sharpness result previously available
only for second-order diffusive terms.
\end{abstract}



\section{Introduction}

Boundary feedback control of one-dimensional hyperbolic systems of conservation and balance laws is a natural framework
for regulating many physical networks in which the state can be accessed only through the endpoints of the spatial domain, see \cite{BastinCoronBook2016} where many examples are studied. Dissipative boundary conditions that guarantee the exponential stability of such systems have been characterized by means of explicit strict Lyapunov functions \cite{CoronBastinAndreaNovel2008,CoronBastin2015}, while boundary feedback laws acting at only one end of the domain have been designed by backstepping \cite{CoronVazquezKrsticBastin2013, CoronNguyen2019}. When the actuator and the sensor share the same
boundary, output-feedback stabilization is by now comparatively well understood
\cite{BastinCoronBook2016,Coron2007,Aam13, AA17b, DSP08, Gug15}. The situation changes considerably when sensing and
actuation are anti-located, that is, placed at opposite endpoints of the domain
\cite{ADM16, KGBS2008, DMHKV16, AnfinsenAamo2017}. The associated input-output map is then no longer passive because of
the transport delay. Simple static feedback laws often fail, and stabilization instead
relies on dynamic compensators that incorporate delayed values of the output, either directly or through a state
predictor. This reliance on delay predictors comes at a price that is well documented in the time-delay systems
literature \cite{MichielsNiculescu2007}. The delay used by the controller must, at least formally, exactly match the
(possibly unknown) transport time of the plant: an arbitrarily small mismatch, coming for instance from an imprecise
characteristic velocity, can move closed-loop poles into the right half-plane and destabilize an otherwise exponentially
stable equilibrium. The sensitivity of exponential stability to the propagation speeds is in fact a classical issue for hyperbolic systems: for linear systems with boundary reflections, Silkowski's criterion characterizes the stability that is robust with respect to the speeds of propagation, see \cite[Section~3]{CoronBastinAndreaNovel2008}. Boundary feedback laws involving time delays have also been analyzed directly, for instance for quasi-linear balance laws with a time-varying boundary delay \cite{Gugat2025}. Moreover, finite-time stability, which is the property behind dead-beat designs, is in general lost when the model is slightly perturbed, although an exponential decay may survive \cite{GugatPerrollazRosier2018}. Whether this fragility is an intrinsic feature of anti-located output feedback, or merely an
artifact of idealized, purely convective plant models, is the structural question addressed in \cite{BCH2022} and
\cite{BCH2025Viscosity}, and it is the starting point of the present paper.
\\

In \cite{BCH2022}, the authors addressed this question for a $2\times2$ example consisting of two transport equations
interconnected in a loop, with the control acting at one boundary and the output measured at the other. The loop is
closed by a dead-beat delayed output feedback of the kind already used for convection-dominated problems
\cite{BreschPietriKrstic2016}. While this feedback is arbitrarily sensitive to velocity mismatches for the inviscid
transport model, they showed that robustness is restored as soon as 
a second-order diffusion term is added to the plant, with an arbitrarily small 
 difussion coefficient $\eta$. Their argument exploits an explicit reduction of the resulting characteristic
equation to a quadratic equation in the spatial eigenvalue $\lambda$, which makes it possible to track the closed-loop
spectrum directly as $\eta\to0^+$. The same paper also showed that this mechanism is not universal, exhibiting a
structurally simpler example in which an arbitrarily small diffusion term destabilizes, rather than stabilizes, the closed
loop. In \cite{BCH2025Viscosity}, the same authors sharpened this picture on a closely related fluid-flow model, a
linearized Saint-Venant type system with the flow rate as boundary control and the fluid density as boundary output, again
closed by a dead-beat delayed feedback and regularized by a second-order viscous term. There, robustness is not only
recovered but quantified exactly: the maximal spectral abscissa of the closed-loop system is shown to converge to the
explicit value $-\ln(2)$ as $\eta\to0^+$, and this exact margin is proved to persist under small perturbations of the
characteristic velocity. The sharpness of this limit, as opposed to a mere robustness statement, is obtained by combining
the explicit quadratic root formula with a topological degree argument, which pins down $-\ln(2)$ simultaneously as an
upper and a lower bound for the limiting spectral abscissa.
\\

In this paper, we pursue this line of research by replacing the second-order diffusion term of
\cite{BCH2022,BCH2025Viscosity} with a fourth-order one. Two difficulties arise immediately. First, a fourth-order spatial
operator requires four boundary conditions per channel instead of two, and choosing them consistently with the hyperbolic
structure of the underlying transport equations already demands some care. Second, and more importantly, the change of
variable that reduces the second-order characteristic equation of \cite{BCH2022,BCH2025Viscosity} to an explicitly
solvable quadratic one no longer applies here: the resulting characteristic equation is a genuine quartic equation in
$\lambda$, whose four roots admit no closed-form expression amenable to the same kind of direct manipulation.
Consequently, the strategy of \cite{BCH2025Viscosity} cannot be replicated as it stands and has to be adapted at the level
of the asymptotic analysis of the roots themselves. Our main contribution is to replace the explicit root formulas of the
second-order case by an asymptotic classification of the four characteristic roots as $\eta\to0^+$, in the spirit of
classical analytic perturbation theory for algebraic equations \cite{Kat66}, precise enough to produce a first-order
expansion of the transfer function together with an explicit, controllable error term. This expansion is then plugged into
a topological degree argument analogous to the one used in \cite{BCH2025Viscosity} to pin down the limiting maximal
spectral abscissa. Our main result is that, despite the lack of usable root formulas, the same limiting value $-\ln(2)$
persists for the fourth-order regularization, and that this sharp margin is itself robust to small perturbations of the
transport velocity. This extends to a genuinely quartic characteristic equation a sharpness result that was previously
available only in the explicitly solvable second-order case, and suggests that the phenomenon uncovered in
\cite{BCH2022,BCH2025Viscosity} is not an artifact of the quadratic structure, but a more robust feature of dissipative
regularizations of anti-located delayed output feedback.
\\

The transfer function approach has a vast literature, since it has been fundamental for understanding linear control
systems in engineering applications. For infinite-dimensional systems described by partial differential equations, we
refer the interested reader to \cite{CurtainMorris2009, Weiss1989,CurtainZwart1995, Staffans2005, TucsnakWeiss2009,TucsnakWeiss2014} and the
references therein.
\\

The remainder of this paper is organized as follows. In Section~\ref{sec: open-loop}, we present the control system and
state our main results concerning transfer function stability and its robustness with respect to errors in the transport
velocity. In Section~\ref{sec: asymptotic}, we develop the asymptotic analysis of the roots of the quartic equation and of
the single-channel transfer function, and prove Theorem~\ref{thm: stabilitysharpness}. This section starts with a heuristic overview of the proof. The estimates are uniform in the
frequency, which requires a separate treatment of moderate and large frequencies, and the sharpness of the bound follows
from a topological degree argument. Section~\ref{sec: robustness} is devoted to proving that the stabilization is robust
with respect to perturbations of the transport velocity. Many proofs are divided into steps to improve readability. Concluding remarks are given in Section~\ref{sec: conclusion}. The proofs of several auxiliary lemmas, which are short computations, are collected in the Appendix.

\section{The control system and main results}\label{sec: open-loop}
 
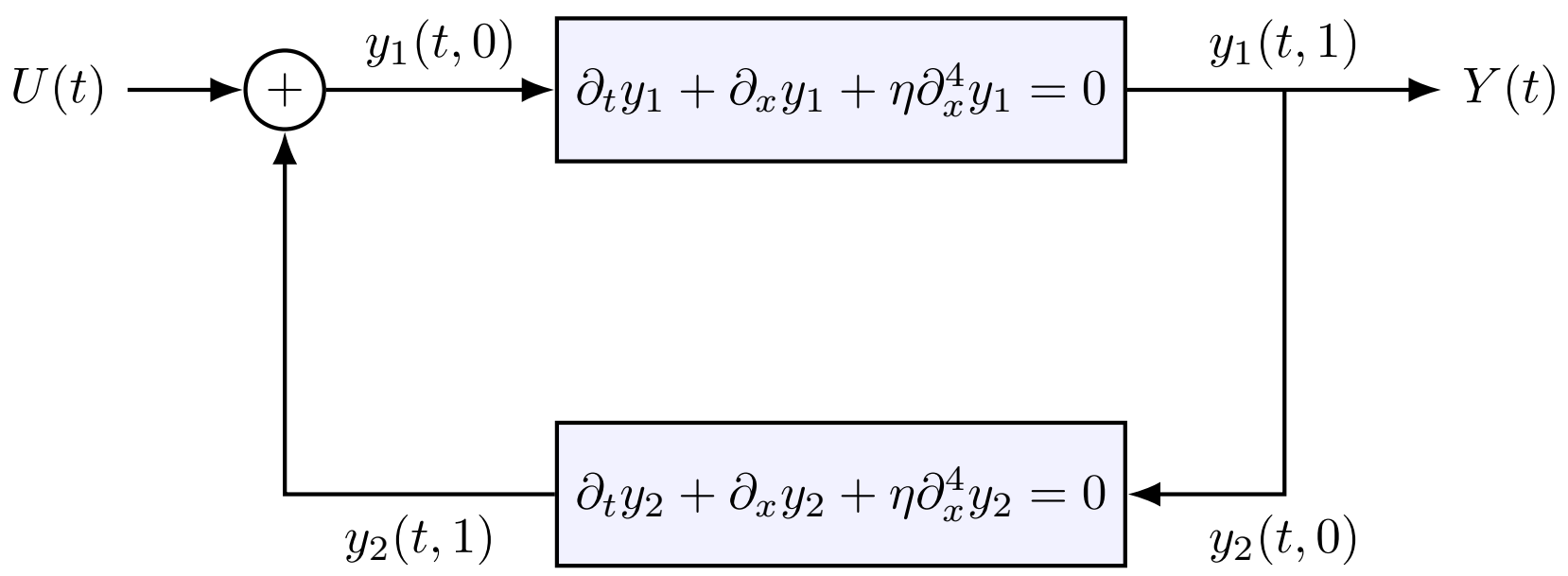
\begin{figure}[h]
\centering
\begin{tikzpicture}[
    >=Latex,
    thick,
    block/.style={draw, minimum height=1cm, minimum width=2.2cm, fill=blue!5, font=\small},
    sum/.style={draw, circle, minimum size=0.55cm, inner sep=0pt, font=\small},
    every node/.style={font=\small}
]
    \node[sum] (sum) at (0,0) {$+$};
    \node[block, right=1.6cm of sum] (b1) {$\partial_t y_1 + \partial_x y_1 +\eta\partial_{x}^4 y_1= 0$};
    \node[block, below=1.8cm of b1] (b2) {$\partial_t y_2 + \partial_x y_2 +\eta\partial_{x}^4 y_2= 0$};
 
    \draw[->] (sum) -- (b1) node[midway, above] {$y_1(t,0)$};
 
    \coordinate (tap) at ($(b1.east)+(1.1,0)$);
    \draw[->] (b1) -- (tap) node[above] {$y_1(t,1)$}  coordinate (tap) -- ++(1.1,0) node[right] {$Y(t)$};
 
    \draw[->] (tap) |- (b2.east) node[pos=0.5, below] {$y_2(t,0)$};
 
    \draw[->] (b2.west) -| (sum.south) node[pos=0.25, below] {$y_2(t,1)$};
 
    \draw[->] (-1.1,0) node[left] {$U(t)$} -- (sum.west);
 
 
\end{tikzpicture}
\caption{Plant given by system \eqref{eq: fourthordersyst}.}
\label{fig: blockdiagram}
\end{figure}

In this paper, we consider two identical transport equations, regularized by fourth-order dissipative terms and coupled
through their boundary conditions. More precisely, for $\eta>0$ we study the system
\begin{equation}\label{eq: fourthordersyst}
  \begin{cases}
    \partial_t y_1(t,x) + \partial_x y_1(t,x) + \eta\partial_x^4 y_1(t,x)=0, \quad &t>0, \ x\in (0,1),\\
    \partial_t y_2(t,x) + \partial_x y_2(t,x) + \eta\partial_x^4 y_2(t,x)=0, \quad &t>0, \ x\in (0,1),\\
    y_1(t,0)=y_2(t,1)+U(t),\quad &t>0,\\
    y_2(t,0)=y_1(t,1),\quad &t>0,\\
    \partial_x^2y_i(t,0)=\partial_x^2y_i(t,1)=\partial_x^3 y_i(t,1)=0, \, i=1,2,\quad &t>0,\\
    Y(t)=y_1(t,1), \quad &t>0.
  \end{cases}
\end{equation}
Here $(y_1,y_2)$ is the state, $U(t)$ denotes the control, and $Y(t)$ the output of the system, see
Figure~\ref{fig: blockdiagram}. The input-output behaviour of this system is described by the transfer function

\begin{equation}\label{eq: IOTransFun}
  G(s)= \frac{Y(s)}{U(s)},
\end{equation}

\noindent where $U(s)$ and $Y(s)$ are the Laplace transforms in time of $U(t)$ and $Y(t)$ respectively. Notice that the plant separates into two identical channels given by 

\begin{equation}\label{single-eq}
  \begin{cases}
    \partial_t y(t,x) + \partial_x y(t,x) + \eta\partial_x^4 y(t,x)=0, \quad &t>0, \ x\in (0,1),\\
    y(t,0)=\bar U(t),\quad &t>0,\\    \partial_x^2y(t,0)=\partial_x^2y(t,1)=\partial_x^3 y(t,1)=0,\quad &t>0,\\
    \bar Y(t)=y(t,1), \quad &t>0,
  \end{cases}
\end{equation}

\noindent where the input is $\bar U(t)$ and the output is $\bar Y(t)$. Our objective is to study the transfer function for each channel \eqref{single-eq}, defined precisely in Definition \ref{c:def-regular} below, and to use the structure of the plant to give an explicit expression for \eqref{eq: IOTransFun}, in order to obtain a transfer function stability result.

\noindent Consistently with the original hyperbolic behavior of the system, we imposed in our systems
three additional boundary conditions:
\[\partial_x^2 y_i(t,0) =0, \quad \partial_x^2 y_i(t,1)=0, \quad \partial_x^3 y_i(t,1)=0.\]
This particular choice is guided both by numerical implementations and by the requirement that the transfer function of a
single channel with fourth-order diffusion converges as $\eta\to 0^+$ to the inviscid one, which is given by
$f_0(s):=e^{-s}$. Note that the latter is the transfer function of the pure transport equation, that is, without the fourth-oder term and the three additional boundary conditions. This convergence is rigorously established later in
Corollary~\ref{c:cor-compact}.
\\

\noindent Taking the Laplace transform in time of \eqref{single-eq} with zero initial condition, we obtain for every $s\in\C$ the boundary value problem
\begin{equation}\label{c:bvp}
s\,y+\left(1+\eps\right)y'+\eta\,y''''=0\ \text{ in }(0,1),\qquad y(0)=y_0,\qquad y''(0)=y''(1)=y'''(1)=0 ,
\end{equation}
where $y_0:=y(s,0)$ is the transform of the input
$\bar U$. We have included a parameter $\eps>-1$ which models an error in the transport velocity, see
\eqref{single-eq-eps} below where the nominal case is $\eps=0$. The following definition is motivated by the approach to transfer functions in \cite{CurtainMorris2009}.

\begin{definition}[Regular points and transfer function]\label{c:def-regular}
Let $\eta>0$ and $\eps>-1$. A point $s\in\C$ is regular for $(\eta,\eps)$ if problem \eqref{c:bvp} 
has exactly one solution for every $y_0$.  By linearity and uniqueness, this solution is $y_0$ times the solution with $y_0=1$. Hence the quotient $y(s,1)/y(s,0)$ does not
depend on $y_0\neq0$, and we define the transfer function at a regular point $s$ by
\[
f_{\eta,\eps}(s):=\frac{y(s,1)}{y(s,0)},\qquad\text{that is, } f_{\eta,\eps}(s)=y(s,1)\text{ for the solution with }y_0=1 .
\]
We write $f_\eta:=f_{\eta,0}$, and we simply say that $s$ is regular for $\eta$ when $\eps=0$.
\end{definition}

\begin{figure}[h]
\centering
\includegraphics[width=0.7\textwidth]{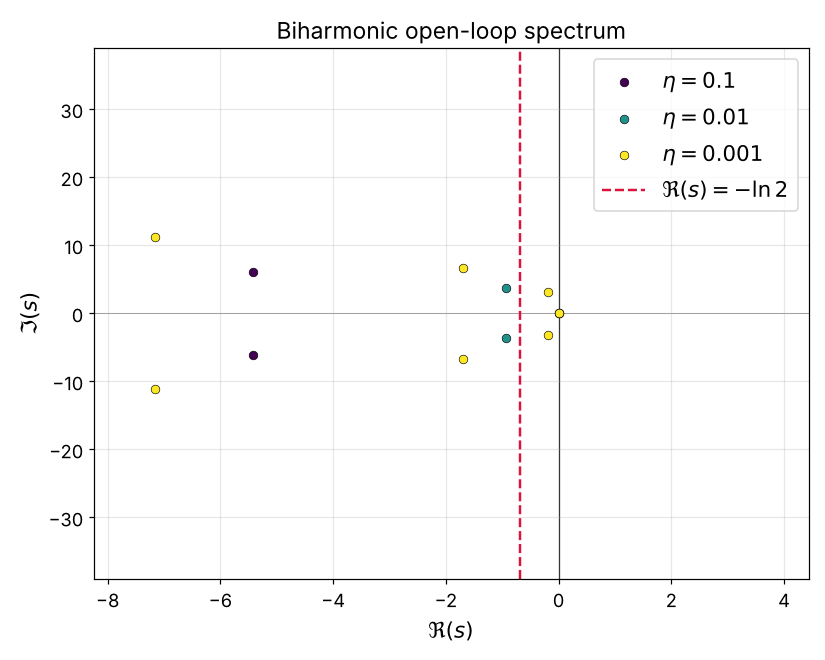}
\caption{Open-loop spectrum of \eqref{eq: fourthordersyst}.}
\label{fig: chap5_openloop}
\end{figure}

\noindent In terms of $f_\eta(s)$, the input-output transfer function of the
open-loop system \eqref{eq: fourthordersyst} is
\begin{equation}\label{eq: overalltransfer}
  G(s) = \frac{Y(s)}{U(s)} = \frac{f_\eta(s)}{1-f_\eta(s)^2}.
\end{equation}
Hence, the poles of the open-loop system are the roots of
\begin{equation}\label{car1}
  f_\eta(s)^2 -1=0.
\end{equation}
Since $f_\eta(0)=1$ for every $\eta>0$, $s=0$ is a pole for every
$\eta>0$, and the open-loop system is not asymptotically stable as can be seen in Figure \ref{fig: chap5_openloop}. This will be shown in Lemma~\ref{c:lem-wp}.

\begin{figure}[h]
\centering
\includegraphics[width=0.7\textwidth]{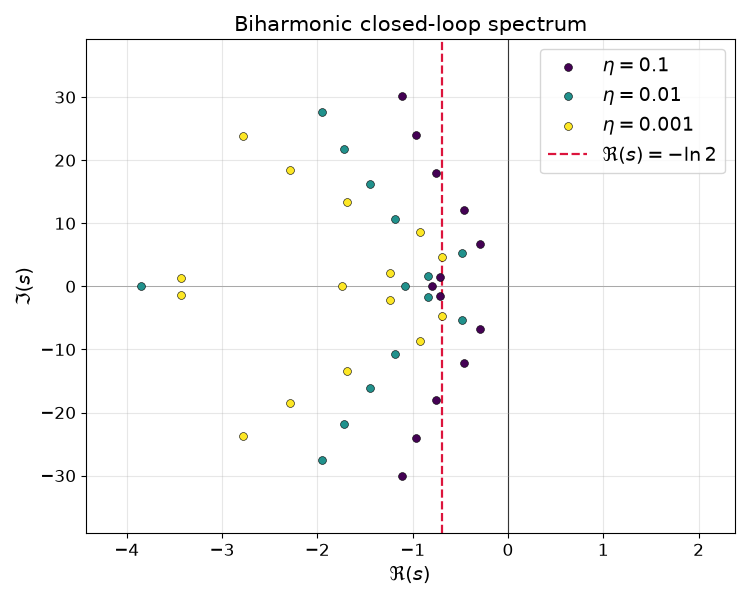}
\caption{Closed-loop spectrum of \eqref{eq: fourthordersyst} with $U(t)=-Y\left(t-1\right)$.}
\label{fig: cerrado}
\end{figure}

\noindent To stabilize this system, we close the loop with the delayed output feedback,
\begin{equation}\label{eq: fourthorderController}
  U(t) = -Y\left(t-1\right).
\end{equation}
The resulting closed-loop spectrum is shown in Figure \ref{fig: cerrado} for some values of $\eta$.

\noindent In terms of $f_\eta(s)$, the closed-loop system reads, in the Laplace domain,
\begin{equation}\label{eq: fourthordersystclosedTransfer}
  \begin{cases}
    y_1(s,1)=f_\eta(s) y_1(s,0),\\
    y_2(s,1)=f_\eta(s) y_2(s,0),\\
    y_1(s,0)=y_2(s,1)-e^{-s}y_1(s,1),\\
    y_2(s,0)=y_1(s,1),\\
    \partial_x^2y_i(s,0)=\partial_x^2y_i(s,1)=\partial_x^3 y_i(s,1)=0, \quad i=1,2,\\
    Y(s)=y_1(s,1).
  \end{cases}
\end{equation}
 
\noindent From the previous input-output transfer function $G(s)$ and the control law \eqref{eq: fourthorderController}, it follows that the closed-loop characteristic equation is
\begin{equation}\label{eq: fourthorderCharEq}
   f_\eta(s)^2 - f_\eta(s)e^{-s} -1=0.
\end{equation}
 
\noindent Our main results concern the location of the roots of the closed-loop characteristic equation
\eqref{eq: fourthorderCharEq}. We say that the closed-loop system is transfer function stable if all these roots lie in a
half-plane $\{\Re s\leq-\alpha\}$ for some $\alpha>0$. This leads to the following definitions.
\begin{definition}\label{def: spectrum}
  The spectrum $\mathcal{S}_\eta$ and maximal spectral abscissa $\sigma_\eta$ of the closed-loop system \eqref{eq: fourthordersystclosedTransfer} are given by
  \[\mathcal{S}_\eta:=\{s\in \C: \ s\text{ is regular for }\eta\text{ and }f_\eta(s)^2 - f_\eta(s)e^{-s} - 1=0\},\]   
  \[\qquad \sigma_\eta:=\sup\{\Re(s): \ s\in \mathcal{S}_\eta\}.\]
\end{definition}
 
\noindent Our first main result is the following. 
 
\begin{theorem}[Transfer function stability and sharpness]\label{thm: stabilitysharpness}
  For every $\delta>0$, there exists $\eta_*>0$ such that
  \begin{equation}\label{eq: fourthordersharpness2}
    -\ln(2) -\delta \leq \sigma_\eta \leq -\ln(2)+\delta, \qquad \forall \eta \in (0, \eta_*).
  \end{equation}
\end{theorem}

\begin{remark}\label{rem: SMT}
  The condition $\sigma_\eta<0$ implies exponential stability of the closed-loop system
\eqref{eq: fourthordersyst}, \eqref{eq: fourthorderController} provided its semigroup satisfies the spectrum-determined
growth condition, e.g.\ through the spectral mapping theorem \cite[page 182]{EN06}. Whether this holds here is an open
problem, as pointed out in \cite{BCH2025Viscosity} for their problem. A possible approach is the Gearhart--Pr\"uss theorem
\cite[Theorem 3.8, page 189]{EN06}.
\end{remark}
 
\noindent Theorem~\ref{thm: stabilitysharpness} establishes the stabilization of the system, as illustrated in Figure~\ref{fig: cerrado}. We now ask whether the feedback law \eqref{eq: fourthorderController} is robust with respect to uncertainties in the transport velocity. To this end, we perturb the nominal velocity $v=1$ by $\varepsilon$ and consider
\begin{equation}\label{eq: fourthordersystepsilon1}
  \begin{cases}
    \partial_t y_1(t,x) + \left(1+\varepsilon\right)\partial_x y_1(t,x) + \eta\partial_x^4 y_1(t,x)=0, \quad &t>0, \ x\in (0,1),\\
    \partial_t y_2(t,x) + \left(1+\varepsilon\right)\partial_x y_2(t,x) + \eta\partial_x^4 y_2(t,x)=0, \quad &t>0, \ x\in (0,1),\\
    y_1(t,0)=y_2(t,1)+U(t),\quad &t>0,\\
    y_2(t,0)=y_1(t,1),\quad &t>0,\\\partial_x^2y_i(t,0)=\partial_x^2y_i(t,1)=\partial_x^3 y_i(t,1)=0, \, i=1,2,\quad &t>0,\\
    Y(t)=y_1(t,1),\quad &t>0.
  \end{cases}
\end{equation}
 
\noindent The perturbed single-channel transfer function $f_{\eta,\varepsilon}(s)$ is defined as the transfer function of the system
\begin{equation}\label{single-eq-eps}
  \begin{cases}
    \partial_t y(t,x) + \left(1+\varepsilon\right)\partial_x y(t,x) + \eta\partial_x^4 y(t,x)=0, \quad &t>0, \ x\in (0,1),\\
    y(t,0)=\hat U(t),\quad &t>0,\\    \partial_x^2y(t,0)=\partial_x^2y(t,1)=\partial_x^3 y(t,1)=0,\quad &t>0,\\
    \hat Y(t)=y(t,1), \quad &t>0,
  \end{cases}
\end{equation}
where the input is $\hat U(t)$ and the output is $\hat Y(t)$.  Its Laplace transform in time is \eqref{c:bvp}, and $f_{\eta,\varepsilon}$ is the function given by Definition~\ref{c:def-regular}. The closed-loop characteristic equation retains the same form:
\begin{equation}\label{eq:char-eq-eps1}
  f_{\eta,\varepsilon}(s)^2 - f_{\eta,\varepsilon}(s)e^{-s} -1=0.
\end{equation}
 
\begin{definition}\label{def: spectrumeps}
  The spectrum $\mathcal{S}_{\eta,\varepsilon}$ and maximal spectral abscissa $\sigma_{\eta,\varepsilon}$ of the closed-loop system \eqref{eq: fourthordersystepsilon1} are given by
  \[\mathcal{S}_{\eta,\varepsilon}:=\{s\in \C: \  s\text{ is regular for }(\eta,\varepsilon)\text{ and }f_{\eta,\varepsilon}(s)^2 - f_{\eta,\varepsilon}(s)e^{-s} - 1=0\},\]
  \[\sigma_{\eta,\varepsilon}:=\sup\{\Re(s): \ s\in \mathcal{S}_{\eta,\varepsilon}\}.\]
\end{definition}
 
\begin{remark}Notice that, by definition, we have 
\begin{equation}
  \mathcal{S}_{\eta,0} = \mathcal{S}_{\eta}, \quad \text{and} \quad \sigma_{\eta,0}= \sigma_{\eta}.
\end{equation}
\end{remark}

\noindent We can now state our second main result.

\begin{theorem}[Robustness]\label{thm:robustness}
For every $\delta>0$ there exists $\eta^*=\eta^*(\delta)>0$ with the following property. For every $\eta\in\left(0,\eta^*\right)$ one
has $\sigma_{\eta,0}\leq-\ln(2)+\delta$, and there exists $\eps_1=\eps_1(\eta,\delta)>0$ such that
\begin{equation}\label{eq: fourthorderrobustness}
\sigma_{\eta,\eps}\leq-\ln(2)+2\delta,\qquad\forall\eps\in\left(-\eps_1,\eps_1\right).
\end{equation}
\end{theorem}
 
\begin{figure}[h]
\centering
\includegraphics[width=0.8\textwidth]{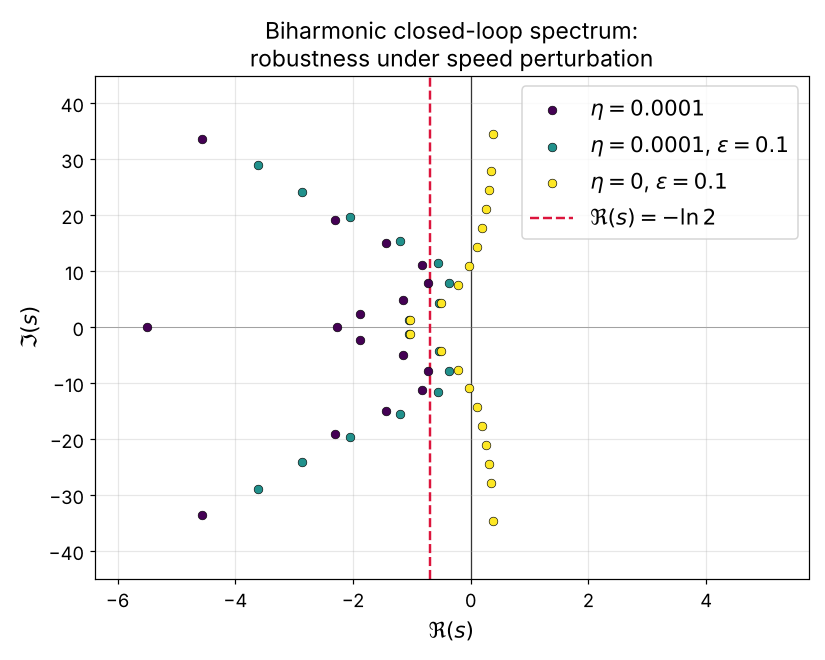}
\caption{Perturbed closed-loop spectrum of \eqref{eq: fourthordersystepsilon1}.}
\label{fig: chap5_perturbedclosedloop}
\end{figure}

\noindent Figure~\ref{fig: chap5_perturbedclosedloop} shows the numerically computed closed-loop spectrum of
\eqref{eq: fourthordersystepsilon1} for a fixed value of $\eta$ and several values of $\varepsilon$, in agreement with
Theorem~\ref{thm:robustness}.
 
\section{Proof of Theorem \ref{thm: stabilitysharpness}}\label{sec: asymptotic}
 
Before going into the details, we explain
heuristically where the value $-\ln(2)$ comes from and how the proof is organized.
 
\medskip\noindent\emph{The inviscid case.} For $\eta=0$, the single-channel transfer function is $e^{-s}$. The input is
transported from $x=0$ to $x=1$ in unit time and the characteristic equation \eqref{eq: fourthorderCharEq} then reads
$e^{-2s}-e^{-2s}-1=0$, that is, $-1=0$. Hence it has no solution. The delay of the controller compensates exactly the transport
time of the plant, and the closed-loop spectrum is empty. This exact compensation is
destroyed by an arbitrarily small error in the transport velocity as shown in \cite{BCH2022}.
 
\medskip\noindent\emph{The fourth-order term.} For $\eta>0$, the solutions of the differential equation in
\eqref{c:bvp} are combinations of exponentials $e^{\lambda x}$, where $\lambda$ solves $\eta\lambda^4+\lambda+s=0$. It is important to note that, even with Ferrari's formula, the four solutions are hard to treat explicitly, so here we employ an asymptotic expansion. Indeed, after expanding the four roots for $\eta$
small, one root stays close to the inviscid one, $\lambda_1(s)\approx-s-\eta s^4$, while the other three are of order
$\eta^{-1/3}$ and produce boundary layers of width of order $\eta^{1/3}$ at $x=0$ and $x=1$. These boundary layers modify the
output only slightly. Therefore, as long as $\eta^{1/3}|s|$ is small,
\begin{equation}\label{c:heur-f}
f_\eta(s)\approx e^{\lambda_1(s)}\approx e^{-s}\,e^{-\eta s^4}.
\end{equation}
For bounded $s$ the factor $e^{-\eta s^4}$ is close to $1$, so that $f_\eta(s)\approx e^{-s}$ and, as in the inviscid case,
\eqref{eq: fourthorderCharEq} has no roots there. On a vertical line $s=\sigma+i\omega$ with $|\omega|$ large, however,
$\eta s^4\approx\eta\omega^4$, and the factor $e^{-\eta s^4}$ is no longer close to $1$ when $t:=\eta\omega^4$ is of order one,
that is, at frequencies $|\omega|$ of order $\eta^{-1/4}$. The fourth-order term acts as a filter which damps the high
frequencies, and this damping is what creates closed-loop poles.
 
\medskip\noindent\emph{Where $-\ln(2)$ comes from.} Let $F:=e^{s}f_\eta(s)$, so that $F\approx e^{-t}$ by \eqref{c:heur-f}.
Multiplying \eqref{eq: fourthorderCharEq} by $e^{2s}$, it becomes $F^2-F=e^{2s}$. Taking moduli, a root $s=\sigma+i\omega$ with
$t=\eta\omega^4$ satisfies, approximately,
\[
e^{2\sigma}\approx\left|e^{-2t}-e^{-t}\right|=e^{-t}\left(1-e^{-t}\right),\qquad\text{that is,}\qquad
\sigma\approx g(t):=-\frac t2+\frac12\ln\left(1-e^{-t}\right).
\]
An elementary computation shows that $g$ attains its maximum on $\left(0,\infty\right)$ at $t=\ln(2)$, and that
$g\left(\ln(2)\right)=-\ln(2)$. Therefore, we expect that $\sigma_\eta\leq-\ln(2)+o(1)$ as $\eta\to0^+$, and that this bound is
attained by roots with $\eta\omega^4\approx\ln(2)$, that is, $|\omega|\approx\left(\ln(2)/\eta\right)^{1/4}$. In other words, the mechanism is essentially the same as in \cite{BCH2022}, but here the error terms coming from the asymptotic analysis must be controlled.
 
\medskip\noindent\emph{The plan.} Three difficulties have to be addressed.
\begin{enumerate}[itemsep=1pt]
\item The approximation \eqref{c:heur-f} requires a description of the four roots of the quartic and a control of the
boundary layers, uniformly in $s$. Since the roots are not explicit, we use the rescaling $\zeta=\eta^{1/3}s$, which turns the
quartic into $\nu^4+\nu+\zeta=0$, and we expand its roots around $\zeta=0$, see Subsection~\ref{c:sec-roots}. Then we compute the
effect of the boundary layers on the output, which gives \eqref{c:heur-f} with an explicit error term, see Subsection~\ref{c:sec-moderate}. The relevant frequencies $|s|\sim\eta^{-1/4}$ satisfy $|\zeta|\sim\eta^{1/12}\to0$, so that
they are covered by these expansions.
\item For $|s|$ of order $\eta^{-1/3}$ or larger, the three terms of the quartic may have the same size, there is no small
parameter, and the expansions break down. In this range we only prove that $f_\eta(s)$ is exponentially small,
see Subsection~\ref{c:sec-large}. This is enough, because $|f_\eta|$ is bounded from below at the points of the spectrum,
see Lemma~\ref{c:lem-lowerF}, so that no root can lie there.
\item The computation above only involves the modulus of $F^2-F$. It explains the upper bound which we prove by in Subsection~\ref{c:sec-upper} by
contradiction using sequences of roots and compactness. However, it does not show that roots
actually exist near $\Re s=-\ln(2)$, which is what the lower bound requires. To prove existence,
in Subsection~\ref{c:sec-lower} we use the topological degree theory, where the mechanism is essentially the same after taking the proper care of the error estimates coming from the asymptotic analysis of the roots and their consequences on the transfer function expression.
\end{enumerate}
In summary, three frequency regimes appear: bounded $s$, where $f_\eta\approx e^{-s}$ and there are no roots; the resonant
regime $|\Im s|\sim\eta^{-1/4}$, which contains the roots that determine $\sigma_\eta$; and $|s|\gtrsim\eta^{-1/3}$, where $f_\eta$
is exponentially small and there are no roots either.
 
 \subsection{First properties of the transfer function}\label{c:sec-prelim}
 
We first collect elementary properties of the transfer function of Definition~\ref{c:def-regular}. They hold for every
$\eps>-1$, and their proofs, which are short computations, are given in Appendix~\ref{app:lem-wp}.
 
\begin{lemma}[Elementary properties of the transfer function]\label{c:lem-wp}
Let $\eta>0$ and $\eps>-1$.
\begin{enumerate}[itemsep=1pt]
\item $s$ is regular for $(\eta,\eps)$ if and only if \eqref{c:bvp} with $y_0=0$ has only the trivial solution $y\equiv0$.
\item $s$ is regular for $(\eta,\eps)$ if and only if $\bar s$ is, and then $f_{\eta,\eps}(\bar s)=\overline{f_{\eta,\eps}(s)}$.
\item $s=0$ is regular for $(\eta,\eps)$ and $f_{\eta,\eps}(0)=1$.
\end{enumerate}
\end{lemma}

\begin{remark}[Symmetry and regularity]\label{c:rem-spectrum}
By the second point on Lemma~\ref{c:lem-wp}, $\Sp_{\eta,\eps}$
is symmetric with respect to the real axis. By third point on Lemma~\ref{c:lem-wp}, $f_\eta(0)=1$, so $s=0$ is a root of \eqref{car1}
for every $\eta>0$. For $\eta$ small, every point of a given half-plane $\{\Re s\geq-M\}$ is regular for $\eta$, see Corollary~\ref{c:cor-polefree}.
\end{remark}
 
\noindent From now on and until the end of this section, $\varepsilon=0$. A key step is to understand the single-channel transfer function $f_\eta(s)$ and to derive its asymptotic expansion as
$\eta\to0^+$. The differential equation in \eqref{c:bvp} is a linear ordinary differential equation in $x$ with constant
coefficients. When its characteristic roots are pairwise distinct, its solutions are of the form
\begin{equation}\label{sol-laplace}
y(s,x)=\sum_{k=1}^4 A_k(s) e^{\lambda_k(s)x}
\end{equation}
where the $A_k(s)$ are coefficients and $\lambda_1(s),\dots,\lambda_4(s)$ are the roots of
\begin{equation}\label{eq: quarticeq}
  \eta\lambda^4 + \lambda + s=0.
\end{equation}
Rather than using explicit formulas for the roots of \eqref{eq: quarticeq}, we describe them asymptotically as $\eta\to0^+$. Imposing the four boundary conditions from system \eqref{single-eq} on $y(s,x)=\sum_{k=1}^4 A_k(s) e^{\lambda_k(s)x}$ produces the linear system $$M(s) (A_1,\dots, A_4)^\top = b,$$ with
\[M(s) = \begin{pmatrix}
  1 & 1 & 1 & 1\\
  \lambda_1^2 & \lambda_2^2 & \lambda_3^2 & \lambda_4^2\\
  \lambda_1^2e^{\lambda_1} & \lambda_2^2e^{\lambda_2} & \lambda_3^2e^{\lambda_3} & \lambda_4^2e^{\lambda_4}\\
  \lambda_1^3e^{\lambda_1} & \lambda_2^3e^{\lambda_2} & \lambda_3^3e^{\lambda_3} & \lambda_4^3e^{\lambda_4}
\end{pmatrix}, \qquad b = \begin{pmatrix}
  y(s,0)\\ 0 \\ 0 \\ 0
\end{pmatrix}.\]
Cramer's rule yields a lengthy but explicit formula for the transfer function
\begin{equation}\label{eq: fourthordertransferfun}
  f_\eta(s)= \frac{y(s,1)}{y(s,0)},
\end{equation}
depending on the four roots $\lambda_k(s)$. This formula is valid as long as the four roots $\lambda_k(s)$ are
pairwise distinct and the matrix $M(s)$ is invertible. Indeed, if the roots are pairwise distinct, the functions $e^{\lambda_kx}$ form a basis of solutions of the
differential equation, so by the first point of Lemma~\ref{c:lem-wp} $s$ is regular if and only if $M(s)$ is invertible, and then
$f_\eta(s)$ is given by \eqref{eq: fourthordertransferfun}.

 
\subsection{The roots of the quartic}\label{c:sec-roots}
 
We set
\begin{equation}\label{c:rescale}
\lambda:=\eta^{-1/3}\nu,\qquad \zeta:=\eta^{1/3}s,\qquad T:=\eta^{-1/3},
\end{equation}
which transforms \eqref{eq: quarticeq} into
\begin{equation}\label{c:Pnu}
P(\nu,\zeta):=\nu^4+\nu+\zeta=0 .
\end{equation}
Equation \eqref{c:Pnu} depends on $(\eta,s)$ only through $\zeta$. Consequently, every quantity of the proof below that is
homogeneous of degree $0$ in $(\lambda_1,\dots,\lambda_4)$ is a function of $\zeta$ alone, and $\eta$ enters separately only through
the exponentials $e^{\pm\lambda_k}=e^{\pm T\nu_k}$. The natural range of the expansions is therefore $|\zeta|$ small, that
is $|s|\ll\eta^{-1/3}$ as we shall see in the next proposition.
\\
 
\noindent Let $\rho_1:=0$, $\rho_2:=-1$, $\rho_3:=e^{i\pi/3}$, $\rho_4:=e^{-i\pi/3}$ be the roots of $P\left(\cdot,0\right)=\nu\left(\nu^3+1\right)$, and
$\Dbar_r:=\{\zeta\in\C:|\zeta|\leq r\}$.
 
\begin{proposition}[Classification of the roots]\label{c:prop-roots}
There exist $\zeta_0>0$, $C_0>0$ and functions $\nu_1,\dots,\nu_4$, holomorphic in a neighbourhood of
$\Dbar_{\zeta_0}$, such that for every $|\zeta|\leq\zeta_0$:
\begin{enumerate}[itemsep=1pt]
\item $\nu_1(\zeta),\dots,\nu_4(\zeta)$ are the roots of \eqref{c:Pnu}. They are simple, and
  $|\nu_k(\zeta)-\rho_k|\leq 1/8$ for $k=1,\dots,4$.
\item $|\nu_1(\zeta)+\zeta+\zeta^4|\leq C_0|\zeta|^7$.
\item $\big|\nu_k(\zeta)-\rho_k-\tfrac{\zeta}{3}-\tfrac{2\rho_k^2}{9}\zeta^2\big|\leq C_0|\zeta|^3$ for $k=2,3,4$.
\end{enumerate}
Consequently, if $\eta>0$ and $|\eta^{1/3}s|\leq\zeta_0$, the roots of \eqref{eq: quarticeq} are $\lambda_k(s)=\eta^{-1/3}\nu_k\left(\eta^{1/3}s\right)$
and
\begin{align}
\lambda_1(s)&=-s-\eta s^4\left(1+r_1(s)\right),\qquad |r_1(s)|\leq C_0\,\eta|s|^3,\label{c:l1}\\
\lambda_k(s)&=\eta^{-1/3}\rho_k+\frac s3+\frac{2\rho_k^2}{9}\,\eta^{1/3}s^2+R_k(s),\qquad
|R_k(s)|\leq C_0\,\eta^{2/3}|s|^3,\quad k=2,3,4.\label{c:lk}
\end{align}
In particular, uniformly for $s$ in a compact subset of $\C$, as $\eta\to0^+$,
\begin{equation}\label{c:compact-roots}
\begin{gathered}
\lambda_1(s)=-s-\eta s^4+O\left(\eta^2\right),\qquad
\mu_k:=\eta^{1/3}\lambda_k(s)-\rho_k=\frac s3\,\eta^{1/3}+\frac{2\rho_k^2}{9}\,s^2\eta^{2/3}+O(\eta),\\
\lambda_k(s)=\eta^{-1/3}\rho_k+\frac s3+O\left(\eta^{1/3}\right),\qquad k=2,3,4.
\end{gathered}
\end{equation}
\end{proposition}
 
\begin{proof}
1. We have $P(\rho_k,0)=0$ and $\partial_\nu P(\rho_k,0)=4\rho_k^3+1$, which equals $1$ for $k=1$ and $-3$ for $k=2,3,4$
since $\rho_k^3=-1$. By the holomorphic implicit function theorem \cite[Section~6.1]{KrantzParks2013} there exist $\zeta_1>0$ and holomorphic functions $\nu_k$
on $\{|\zeta|<\zeta_1\}$ with $\nu_k(0)=\rho_k$ and $P\left(\nu_k(\zeta),\zeta\right)=0$. By continuity we may fix
$\zeta_0\in\left(0,\zeta_1/2\right)$ such that $|\nu_k(\zeta)-\rho_k|\leq1/8$ for $|\zeta|\leq2\zeta_0$ and all $k$. Since
$|\rho_j-\rho_k|\geq1$ for $j\ne k$, the discs of radius $1/8$ around the $\rho_k$ are pairwise disjoint, thus
$\nu_1(\zeta),\dots,\nu_4(\zeta)$ are four distinct roots of the quartic $P\left(\cdot,\zeta\right)$, hence all its roots and they are simple.
\\

\noindent 2. Writing $\nu_1=-\zeta u$, the equation $P(\nu_1,\zeta)=0$ becomes $u=1+\zeta^3u^4$. The map
$(u,w)\mapsto u-1-wu^4$ vanishes at $(1,0)$ and its derivative with respect to $u$ equals $1$ there. By the holomorphic implicit
function theorem  \cite[Section~6.1]{KrantzParks2013} there is a function $U$, holomorphic near $0$, with $U(0)=1$ and $U(w)=1+wU(w)^4$. Differentiating at
$w=0$ gives $U'(0)=1$. For $\zeta$ small, $-\zeta U\left(\zeta^3\right)$ is a root of $P\left(\cdot,\zeta\right)$ in the disc of radius $1/8$
around $\rho_1=0$, so it coincides with $\nu_1(\zeta)$ by item 1. Since $U(w)-1-w$ vanishes to second order at $w=0$, there
is $C>0$ with $|U(w)-1-w|\leq C|w|^2$ for $w$ small, and therefore
\[
|\nu_1(\zeta)+\zeta+\zeta^4|=|\zeta|\,|U(\zeta^3)-1-\zeta^3|\leq C|\zeta|^7 .
\]
Decreasing $\zeta_0$ if necessary, this holds for $|\zeta|\leq\zeta_0$.
\\

\noindent 3. Let $k\in\{2,3,4\}$. The function $\nu_k$ is holomorphic in a neighbourhood of $\Dbar_{2\zeta_0}$ and satisfies
$|\nu_k|\leq|\rho_k|+\frac18=\frac98$ there. Differentiating $P\left(\nu_k(\zeta),\zeta\right)=0$ twice gives
\[
\left(4\nu_k^3+1\right)\,\nu_k'+1=0,\qquad \left(4\nu_k^3+1\right)\,\nu_k''+12\,\nu_k^2\,(\nu_k')^2=0 .
\]
At $\zeta=0$ we have $4\rho_k^3+1=-3$, hence $\nu_k'(0)=\frac13$ and $\nu_k''(0)=\frac{4\rho_k^2}{9}$, so the Taylor polynomial
of degree $2$ of $\nu_k$ at $0$ is $\rho_k+\frac\zeta3+\frac{2\rho_k^2}{9}\zeta^2$. By the Cauchy estimates on the circle
$|\zeta|=2\zeta_0$, the remainder satisfies, for $|\zeta|\leq\zeta_0$,
\[
\Big|\nu_k(\zeta)-\rho_k-\frac\zeta3-\frac{2\rho_k^2}{9}\zeta^2\Big|\leq\frac{9\left(|\zeta|/(2\zeta_0)\right)^3}{8(1-|\zeta|/(2\zeta_0))}\leq\frac{9}{32\,\zeta_0^3}\,|\zeta|^3 .
\]
We take for $C_0$ the larger of the constants in items 2 and 3.
\\

\noindent Consequently, if $|\eta^{1/3}s|\leq\zeta_0$, the map $\nu\mapsto T\nu$ sends the roots of \eqref{c:Pnu} with
$\zeta=\eta^{1/3}s$ onto those of \eqref{eq: quarticeq}. Since $T\zeta=s$, $T\zeta^2=\eta^{1/3}s^2$, $T\zeta^3=\eta^{2/3}s^3$,
$T\zeta^4=\eta s^4$ and $T\zeta^7=\eta^2s^7$, multiplying item 2 by $T$ gives $\lambda_1=-s-\eta s^4+R$ with
$|R|\leq C_0\eta^2|s|^7$, which is \eqref{c:l1} with $r_1:=-R/\left(\eta s^4\right)$ and $r_1(0):=0$. Multiplying item 3 by $T$ gives
\eqref{c:lk}. If $s$ ranges in a compact set, then $|\eta^{1/3}s|\leq\zeta_0$ for $\eta$ small, and \eqref{c:compact-roots}
follows from \eqref{c:l1}, \eqref{c:lk} and, for $\mu_k=\nu_k(\zeta)-\rho_k$, directly from item 3.
\end{proof}

\begin{remark}
In what follows, we call $\lambda_1$ the finite branch and $\lambda_2,\lambda_3,\lambda_4$ the exploding branches.
\end{remark}

\begin{remark}[Boundary layers]\label{c:rem-layers}
By item 1 of Proposition~\ref{c:prop-roots}, if $|\eta^{1/3}s|\leq\zeta_0$ then $\Re\lambda_2\leq-\frac78\eta^{-1/3}$ and
$\Re\lambda_3, \Re\lambda_4\geq\frac38\eta^{-1/3}$. Thus $e^{\lambda_2x}$ is a boundary layer at $x=0$, while $\lambda_3, \lambda_4$
produce a boundary layer at $x=1$ carried by two modes. This is consistent with the boundary conditions of
\eqref{eq: fourthordersyst} since we have one extra condition at $x=0$ and two at $x=1$.
\end{remark}
 
\noindent The next lemma follows from Vieta's formulas, since \eqref{eq: quarticeq} has no cubic or quadratic term. It will allow
us to keep track of the error terms in the expansion of $f_\eta(s)$.
 
\begin{lemma}[Vieta identity]\label{lem:vieta}
Let $\lambda_1, \lambda_2, \lambda_3, \lambda_4$ be the four roots of \eqref{eq: quarticeq}. Then
\begin{equation}\label{eq:vieta}
\lambda_3+\lambda_4 = -\left(\lambda_1+\lambda_2\right), \qquad\text{and}\qquad \lambda_3^2+\lambda_3\lambda_4+\lambda_4^2 = \lambda_1\lambda_2 .
\end{equation}
\end{lemma}

\noindent Since \eqref{c:Pnu} is obtained from \eqref{eq: quarticeq} by the rescaling \eqref{c:rescale}, its roots
$\nu_k=\eta^{1/3}\lambda_k$ satisfy the same identities.
 
\subsection{The transfer function at moderate frequencies}\label{c:sec-moderate}
 
The functions below collect the effect of the two boundary layers on the output. When the slow mode
$e^{\lambda_1x}$ reaches $x=1$, the conditions $y''(1)=y'''(1)=0$ excite the two growing modes $e^{\lambda_3\left(x-1\right)}$ and
$e^{\lambda_4\left(x-1\right)}$, and these modify the output $y(1)$ by the factor $1-X_1$ (see \eqref{c:Phi} below). The term $X_2$ (see \eqref{c:Phi} again) plays
the same role for the mode $e^{\lambda_2x}$, but it multiplies an exponentially small quantity. At $x=0$, the condition
$y''(0)=0$ excites the decaying mode $e^{\lambda_2x}$, which changes the amplitude of the slow mode by the factor
$1/\left(1-\left(\nu_1/\nu_2\right)^2\right)$. The function $\Phi$ combines both corrections, so that $f_\eta\approx e^{\lambda_1}\Phi$ as we will show in Proposition~\ref{c:prop-A}.
\\
 
\noindent For $|\zeta|\leq\zeta_0$ we define, with $\nu_k=\nu_k(\zeta)$,
\begin{equation}\label{c:Phi}
X_1(\zeta):=\frac{\nu_1^3\left(\nu_1+2\nu_2\right)}{(\nu_3\nu_4)^2},\qquad
X_2(\zeta):=\frac{\nu_2^3\left(\nu_2+2\nu_1\right)}{(\nu_3\nu_4)^2},\qquad
\Phi(\zeta):=\frac{1-X_1(\zeta)}{1-\left(\nu_1(\zeta)/\nu_2(\zeta)\right)^2}.
\end{equation}
Since numerator and denominator of $X_1$ and $X_2$ are homogeneous polynomials of degree $4$ in the roots, the
factor $\eta^{1/3}$ in $\nu_k=\eta^{1/3}\lambda_k$ cancels, and the same formulas hold with the roots of \eqref{eq: quarticeq}:
$X_1=\lambda_1^3\left(\lambda_1+2\lambda_2\right)/(\lambda_3\lambda_4)^2$ and $X_2=\lambda_2^3\left(\lambda_2+2\lambda_1\right)/(\lambda_3\lambda_4)^2$. In
particular, $X_1$, $X_2$ and $\Phi$ depend on $(\eta,s)$ only through $\zeta$.
 
\begin{proposition}[Uniform expansion of $f_\eta$ at moderate frequencies]\label{c:prop-A}
Let $\zeta_0$ be as in Proposition~\ref{c:prop-roots}. The functions $X_1,X_2,\Phi$ are holomorphic in a neighbourhood of
$\Dbar_{\zeta_0}$, and there is $C_\Phi>0$ such that, for $|\zeta|\leq\zeta_0$,
\begin{equation}\label{c:Phi-bounds}
|\Phi(\zeta)-1-\zeta^2|\leq C_\Phi|\zeta|^3,\qquad \tfrac12\leq|\Phi(\zeta)|\leq2,\qquad
|X_1(\zeta)|\leq C_\Phi|\zeta|^3,\qquad X_2(0)=1 .
\end{equation}
Moreover, there exist $\eta_A\in(0,1]$ and $C_A>0$ such that, for every $\eta\in(0,\eta_A]$ and every $s\in\C$ with
$|\eta^{1/3}s|\leq\zeta_0$, the matrix $M(s)$ is invertible, so $s$ is regular for $\eta$, and
\begin{equation}\label{c:lemmaA}
\begin{gathered}
f_\eta(s)=e^{\lambda_1(s)}\,\Phi\left(\eta^{1/3}s\right)\,\left(1+r_2(s)\right)+r_3(s),\\
|r_2(s)|\leq C_Ae^{-\eta^{-1/3}/4}\leq\tfrac12,\qquad |r_3(s)|\leq C_Ae^{-\eta^{-1/3}/2}.
\end{gathered}
\end{equation}
Moreover, the map $(\eta,s)\mapsto f_\eta(s)$ is continuous on the set
$\{(\eta,s):\ 0<\eta\leq\eta_A,\ |\eta^{1/3}s|\leq\zeta_0\}$.
\end{proposition}
 
\begin{proof}
Recall $T=\eta^{-1/3}$, $\zeta=\eta^{1/3}s$, $\lambda_k=T\nu_k(\zeta)$, and $C$ denotes constants depending only on
$\zeta_0$.
\\
 
\noindent\emph{Step 0: uniform bounds.} By Proposition~\ref{c:prop-roots}, for $|\zeta|\leq\zeta_0$:
\begin{equation}\label{c:step0}
\begin{gathered}
|\nu_1|\leq\tfrac18,\quad \Re\nu_2\leq-\tfrac78,\quad \Re\nu_3,\ \Re\nu_4\geq\tfrac38,\quad
\Re\left(\nu_i-\nu_j\right)\leq-\tfrac14\ \ (i\in\{1,2\},\ j\in\{3,4\}),\\
\tfrac78\leq|\nu_k|\leq\tfrac98\ \left(k=2,3,4\right),\quad |\nu_3-\nu_4|\geq\sqrt3-\tfrac14\geq1,\quad
|\nu_2^2-\nu_1^2|\geq\tfrac34,\quad |\nu_1/\nu_2|\leq\tfrac17 .
\end{gathered}
\end{equation}
 
\noindent\emph{Step 1: the layer at $x=1$.} Since the roots are simple, we write the solution with $y(0)=1$ as
\[
y(x)=A_1e^{\lambda_1x}+A_2e^{\lambda_2x}+B_3e^{\lambda_3\left(x-1\right)}+B_4e^{\lambda_4\left(x-1\right)},
\]
where the last two coefficients are normalized at $x=1$. With $c_i:=A_ie^{\lambda_i}$, the conditions $y''(1)=y'''(1)=0$ read
$\lambda_3^2B_3+\lambda_4^2B_4=-\sum_ic_i\lambda_i^2$ and $\lambda_3^3B_3+\lambda_4^3B_4=-\sum_ic_i\lambda_i^3$. This $2\times2$
system has determinant $\lambda_3^2\lambda_4^2\left(\lambda_4-\lambda_3\right)=T^5\nu_3^2\nu_4^2\left(\nu_4-\nu_3\right)\neq0$, and Cramer's rule gives,
after dividing numerators and denominators by $T^3$,
\begin{equation}\label{c:B34}
\begin{gathered}
B_3=\sum_{i=1,2}c_i\beta_{3i},\qquad B_4=\sum_{i=1,2}c_i\beta_{4i},\\
\beta_{3i}:=\frac{\nu_i^2\left(\nu_4-\nu_i\right)}{\nu_3^2\left(\nu_3-\nu_4\right)},\qquad
\beta_{4i}:=-\frac{\nu_i^2\left(\nu_3-\nu_i\right)}{\nu_4^2\left(\nu_3-\nu_4\right)} .
\end{gathered}
\end{equation}
By \eqref{c:step0}, $|\beta_{ji}|\leq4$. Moreover,
\[
\beta_{3i}+\beta_{4i}=\frac{\nu_i^2\left[\nu_4^2\left(\nu_4-\nu_i\right)-\nu_3^2\left(\nu_3-\nu_i\right)\right]}{(\nu_3\nu_4)^2\left(\nu_3-\nu_4\right)}
=-\frac{\nu_i^2\left[\nu_3^2+\nu_3\nu_4+\nu_4^2-\nu_i\left(\nu_3+\nu_4\right)\right]}{(\nu_3\nu_4)^2},
\]
and by Lemma~\ref{lem:vieta} the bracket equals $\nu_1\nu_2+\nu_i\left(\nu_1+\nu_2\right)$, that is $\nu_1\left(\nu_1+2\nu_2\right)$ for $i=1$ and
$\nu_2\left(\nu_2+2\nu_1\right)$ for $i=2$. The large roots $\nu_3, \nu_4$ have disappeared from the bracket, and
$\beta_{3i}+\beta_{4i}=-X_i$, that is,
\begin{equation}\label{c:B3plusB4}
B_3+B_4=-c_1X_1-c_2X_2 .
\end{equation}
The functions $X_1, X_2$ are holomorphic near $\Dbar_{\zeta_0}$, since their denominators do not vanish by \eqref{c:step0},
$\nu_1^3$ has a zero of order $3$ at $\zeta=0$, so $|X_1(\zeta)|\leq C|\zeta|^3$ and $X_2(0)=1$. In the
notation of the original variables, $X_1=O\left(\eta|s|^3\right)$ and $X_2=1+O\left(\eta^{1/3}|s|\right)$. In particular $X_2$ is not small,
but it only multiplies $A_2e^{\lambda_2}$, which is exponentially small.
\\
 
\noindent\emph{Step 2: the layer at $x=0$ as a perturbed $2\times2$ system.} By \eqref{c:B34},
$B_je^{-\lambda_j}=\sum_iA_i\beta_{ji}e^{T\left(\nu_i-\nu_j\right)}$ for $j=3, 4$. Substituting into $y(0)=1$ and $y''(0)=0$, and
dividing the latter by $T^2$, we obtain
\begin{equation}\label{c:2x2}
\left(1+\eps_{11}\right)A_1+\left(1+\eps_{12}\right)A_2=1,\qquad \left(\nu_1^2+\eps_{21}\right)A_1+\left(\nu_2^2+\eps_{22}\right)A_2=0,
\end{equation}
with
\[
\eps_{1i}:=\beta_{3i}e^{T\left(\nu_i-\nu_3\right)}+\beta_{4i}e^{T\left(\nu_i-\nu_4\right)},\qquad
\eps_{2i}:=\nu_3^2\beta_{3i}e^{T\left(\nu_i-\nu_3\right)}+\nu_4^2\beta_{4i}e^{T\left(\nu_i-\nu_4\right)} .
\]
By \eqref{c:step0}, $|e^{T\left(\nu_i-\nu_j\right)}|\leq e^{-T/4}$, hence $|\eps_{ji}|\leq Ce^{-T/4}$. The determinant of \eqref{c:2x2}
is $\Delta=\left(\nu_2^2-\nu_1^2\right)+\vartheta$, where the lower-order term $\vartheta$ collects the products containing some $\eps_{ji}$ and satisfies $|\vartheta|\leq Ce^{-T/4}$. Choose $\eta_A$ so small that
$Ce^{-T/4}\leq\frac14$ for $\eta\leq\eta_A$. Then $|\Delta|\geq\frac12$ and \eqref{c:2x2} is uniquely solvable:
\[
A_1=\frac{\nu_2^2+\eps_{22}}{\Delta}=\frac{1}{1-\left(\nu_1/\nu_2\right)^2}\,\left(1+r_2\right),\qquad
1+r_2:=\frac{1+\eps_{22}/\nu_2^2}{1+\vartheta/\left(\nu_2^2-\nu_1^2\right)},\qquad
A_2=-\frac{\nu_1^2+\eps_{21}}{\Delta},
\]
so that $|r_2|\leq Ce^{-T/4}$ and $|A_2|\leq C$.
\\
 
\noindent\emph{Step 3: Invertibility of $M(s)$.} If $y(0)=0$, Step 1 still expresses $(B_3,B_4)$ through $(A_1,A_2)$, and Step 2 gives the
system \eqref{c:2x2} with right-hand side $0$, whence $A_1=A_2=0$ and $B_3=B_4=0$. Thus the homogeneous problem has
only the trivial solution. Together with Step~2 for $y(0)=1$ and linearity, \eqref{c:bvp} has exactly one
solution for every $y_0$. Hence $s$ is regular for $\eta$ and $M(s)$ is invertible.
\\
 
\noindent\emph{Step 4: the output.} By \eqref{c:B3plusB4},
\[
f_\eta(s)=y(1)=A_1e^{\lambda_1}+A_2e^{\lambda_2}+B_3+B_4=A_1e^{\lambda_1}\left(1-X_1\right)+A_2e^{\lambda_2}\left(1-X_2\right).
\]
By Step~2 this is \eqref{c:lemmaA} with $r_3:=A_2e^{\lambda_2}\left(1-X_2\right)$. Since $|e^{\lambda_2}|=e^{T\Re\nu_2}\leq e^{-7T/8}$,
we get $|r_3|\leq Ce^{-T/2}$. Enlarging $C$ and decreasing $\eta_A$ we obtain the constant $C_A$ and the bound
$|r_2|\leq\frac12$ in \eqref{c:lemmaA}.
\\
 
\noindent\emph{Step 5: estimates for $\Phi$.} By \eqref{c:step0}, $|\nu_1/\nu_2|<1$, so the denominator of $\Phi$
does not vanish and $\Phi$ is holomorphic near $\Dbar_{\zeta_0}$. Since $\nu_1(0)=0$ and $X_1(0)=0$, we have $\Phi(0)=1$.
Decreasing $\zeta_0$ if necessary, $|\Phi(\zeta)-1|\leq\frac12$ on $\Dbar_{\zeta_0}$, and hence $\frac12\leq|\Phi|\leq2$. By Proposition~\ref{c:prop-roots},
$\nu_1=-\zeta+O\left(\zeta^4\right)$ and $\nu_2=-1+\frac\zeta3+O\left(\zeta^2\right)$, so $\nu_1/\nu_2=\zeta+O\left(\zeta^2\right)$ and
$\left(1-\left(\nu_1/\nu_2\right)^2\right)^{-1}=1+\zeta^2+O\left(\zeta^3\right)$. Together with $X_1=O\left(\zeta^3\right)$ this gives $\Phi(\zeta)=1+\zeta^2+O\left(\zeta^3\right)$
as $\zeta\to0$. Hence $\left(\Phi(\zeta)-1-\zeta^2\right)/\zeta^3$ has a removable singularity at $0$ and is bounded on
$\Dbar_{\zeta_0}$, which is the first bound in \eqref{c:Phi-bounds}.
\\
 
\noindent\emph{Step 6: continuity.} On the set $\{(\eta,s):\ 0<\eta\leq\eta_A,\ |\eta^{1/3}s|\leq\zeta_0\}$, the roots
$\nu_k\left(\eta^{1/3}s\right)$, the number $T=\eta^{-1/3}$, the coefficients $\beta_{ji}$ of \eqref{c:B34} and the coefficients
$\eps_{ji}$ of \eqref{c:2x2} are continuous functions of $(\eta,s)$, and the determinant $\Delta$ of \eqref{c:2x2} satisfies
$|\Delta|\geq\frac12$. By Cramer's rule, $A_1$ and $A_2$ are continuous, hence so are $B_3$ and $B_4$ by \eqref{c:B34}, and so is
$f_\eta(s)=A_1e^{\lambda_1}+A_2e^{\lambda_2}+B_3+B_4$.
\end{proof}

\noindent Proposition \ref{c:prop-A} has the following two consequences for the compact case and resonant scale, whose proofs could be referred in Appendices \ref{app:cor-compact} and \ref{app:cor-resonant}, respectively.
 
\begin{corollary}[$s$ in a compact set]\label{c:cor-compact}
Let $K\subset\C$ be compact. Then, uniformly for $s\in K$, as $\eta\to0^+$,
\begin{equation}\label{c:fcompact1}
f_\eta(s)=e^{\lambda_1(s)}\left(1+\eta^{2/3}s^2+O(\eta)\right)=e^{\lambda_1(s)}\left(1+O\left(\eta^{2/3}\right)\right),
\end{equation}
\begin{equation}\label{c:fcompact2}
  f_\eta(s)=e^{-s}\left(1+O\left(\eta^{2/3}\right)\right).
\end{equation}
In particular $f_\eta\to e^{-s}$ uniformly on $K$.
\end{corollary}

\begin{corollary}[Resonant scale, uniformly and with a rate]\label{c:cor-resonant}
Let $Q\subset\R\times\left(0,\infty\right)$ be compact. There exist $\eta_Q\in(0,\eta_A]$ and $C_Q>0$ such that, for every
$\eta\in(0,\eta_Q]$ and every $(\sigma,\psi)\in Q$, the point $s:=\sigma+i\eta^{-1/4}\psi$ satisfies
$|\eta^{1/3}s|\leq\zeta_0$, so $s$ is regular for $\eta$ by Proposition~\ref{c:prop-A}, and
\begin{equation}\label{c:resonant}
|\eta s^4-\psi^4|\leq C_Q\,\eta^{1/4},\qquad \big|e^{s}f_\eta(s)-e^{-\psi^4}\big|\leq C_Q\,\eta^{1/6}.
\end{equation}
In particular $f_\eta(s)=e^{-s-\psi^4}\left(1+O\left(\eta^{1/6}\right)\right)$ uniformly on $Q$.
\end{corollary}

\subsection{The transfer function at large frequencies}\label{c:sec-large}
 
Proposition~\ref{c:prop-A} controls $f_\eta$ only for $|\eta^{1/3}s|\leq\zeta_0$. The upper bound of
Theorem~\ref{thm: stabilitysharpness} requires a control in the whole half-plane $\{\Re s\geq-M\}$ for $M>0$, which contains the
frequencies $|s|\gtrsim\eta^{-1/3}$ and the case $\Re s\to+\infty$. In this range there is no small parameter: when
$|\eta^{1/3}s|$ is of order one, the three terms of $\eta\lambda^4+\lambda+s$ have the same size, and the classification of
Proposition~\ref{c:prop-roots} is no longer available. The result below replaces the asymptotic expansions by a
compactness argument. After the rescaling $\lambda=N\xi$, $N=\left(|s|/\eta\right)^{1/4}$, the quartic becomes
$\xi^4+\kappa\xi+e^{i\theta}=0$, whose coefficients $(\theta,\kappa)$ range in a compact set. Three qualitative facts then yield uniform constants: no
root on the imaginary axis, two roots on each side of $\{\Re \xi =0\}$ and the two roots on the right are simple. Therefore the two boundary
layers decouple exactly as in Proposition~\ref{c:prop-A}. The conclusion is simply that $|f_\eta|$ is exponentially small
in $N$, which is enough, since $|f_\eta|$ is bounded from below at points of the spectrum as we will see in Lemma~\ref{c:lem-lowerF} below.
\\
 
\begin{proposition}[Large frequencies]\label{c:prop-B}
For $\eta>0$ and $s\in\C\setminus\{0\}$ set
\begin{equation}\label{c:Nkappa}
\theta:=\arg s\in(-\pi,\pi],\qquad N:=\left(\frac{|s|}{\eta}\right)^{1/4},\qquad {\kappa:=|\eta^{1/3}s|^{-3/4}.}
\end{equation}
Let $\kappa_0>0$. There exist $\theta_1\in\left(0,\pi/6\right)$, $N_0\geq1$ and $c, C>0$, depending only on $\kappa_0$, such that:
if $|\theta|\leq\frac\pi2+\theta_1$, $\kappa\leq\kappa_0$ and $N\geq N_0$, then problem \eqref{c:bvp} is uniquely solvable for every $y_0$, i.e.\ $s$ is regular for $\eta$, and its solution with $y_0=1$ satisfies $|y(1)|\leq Ce^{-cN}$. In particular
$|f_\eta(s)|\leq Ce^{-cN}$.
\end{proposition}
 
\begin{proof}

\noindent\emph{Step 1: rescaled quartic.} Substituting $\lambda=N\xi$ in $\eta\lambda^4+\lambda+s=0$ and dividing by
$\eta N^4=|s|$ we get
\begin{equation}\label{c:Qxi}
Q_{\theta,\kappa}(\xi):=\xi^4+\kappa\,\xi+e^{i\theta}=0,
\end{equation}
since $N/|s|=1/\left(\eta N^3\right)=\kappa$. If $\xi_1,\dots,\xi_4$ are the roots of \eqref{c:Qxi}, then $\lambda_k:=N\xi_k$ are the
roots of $\eta\lambda^4+\lambda+s=0$. In particular $\Re\lambda_k=N\Re\xi_k$.
\\
 
\noindent\emph{Step 2: the parameter set.} Fix $\theta_1\in\left(0,\pi/6\right)$ with $\cos\theta_1>\kappa_0\left(\sin\theta_1\right)^{1/4}$ and let
$K:=\left[-\frac\pi2-\theta_1,\frac\pi2+\theta_1\right]\times[0,\kappa_0]$, a compact connected set.
 
\begin{enumerate}
  \item \emph{No root of $Q_{\theta,\kappa}$ lies on $i\R$ for $(\theta,\kappa)\in K$.} If $\xi=iy$ with $y\in\R$, the real and
imaginary parts of $Q_{\theta,\kappa}(iy)=y^4+i\kappa y+e^{i\theta}=0$ give $y^4=-\cos\theta$ and $\kappa y=-\sin\theta$.
If $|\theta|<\frac\pi2$ the first identity is impossible. If $|\theta|=\frac\pi2$ it gives $y=0$ and then $\sin\theta=0$,
which is impossible. If $\frac\pi2<|\theta|\leq\frac\pi2+\theta_1$, then $|y|=|\cos\theta|^{1/4}$ and
$\kappa=|\sin\theta|/|\cos\theta|^{1/4}\geq\cos\theta_1/\left(\sin\theta_1\right)^{1/4}>\kappa_0$, a contradiction.

  \item \emph{Two roots, $\xi_1,\xi_2$, lie in $\{\Re\xi<0\}$ and two, $\xi_3,\xi_4$, in $\{\Re\xi>0\}$.} The roots of a monic
  polynomial depend continuously on its coefficients, as an unordered $4$-tuple by Rouch\'e's theorem. Since
  $K$ is connected and, by item 1, no root crosses $i\R$, the number of roots in $\{\Re\xi<0\}$ is constant on $K$. At
  $(\theta,\kappa)=(0,0)$ the roots are $e^{\pm i\pi/4}$ and $e^{\pm3i\pi/4}$.

  \item \emph{$\xi_3\neq\xi_4$.} A double root satisfies $Q_{\theta,\kappa}(\xi)=Q'_{\theta,\kappa}(\xi)=0$. If $\kappa=0$,
  then $Q'=4\xi^3=0$ forces $\xi=0$, but $Q(0)=e^{i\theta}\neq0$. If $\kappa>0$, then $4\xi^3=-\kappa$ gives
  $\xi^4=-\kappa\xi/4$, so $0=Q(\xi)=\frac34\kappa\xi+e^{i\theta}$, i.e.\ $\xi=-\frac{4}{3\kappa}e^{i\theta}$, and then
  $4\xi^3=-\kappa$ yields $e^{3i\theta}=\frac{27\kappa^4}{256}>0$, i.e.\ $\theta\in\{0,\pm\frac{2\pi}3\}$ and, since
  $|e^{3i\theta}|=1$, $\kappa^4=\frac{256}{27}$. The values
  $\pm\frac{2\pi}3$ are not in $K$ because $\theta_1<\pi/6$. For $\theta=0$ the double root $\xi=-\frac4{3\kappa}$ is real
  and negative. Hence a double root can only occur in the left pair, namely for $\theta=0$ and when $\kappa$ attains $\kappa_c= 4 \cdot 3^{-3/4}$.

  \item \emph{Uniform constants.} The functions $(\theta,\kappa)\mapsto\min_k|\Re\xi_k|$ and
  $(\theta,\kappa)\mapsto|\xi_3-\xi_4|$, where the right pair is well defined as an unordered pair by point 2, are continuous on the
  compact set $K$, and positive by points 1 and 3. Let $c_K>0$ and $d_K>0$ be lower bounds. By Cauchy's bound,
  $|\xi_k|\leq1+\max\{\kappa,1\}\leq C_K:=2+\kappa_0$. Consequently, for $(\theta,\kappa)\in K$,
  \begin{equation}\label{c:unifxi}
  \begin{gathered}
    \Re\xi_{1,2}\leq-c_K,\qquad \Re\xi_{3,4}\geq c_K,\qquad |\xi_k|\leq C_K,\\
    |\xi_3|,|\xi_4|\geq c_K,\qquad |\xi_3-\xi_4|\geq d_K,\qquad |\xi_1+\xi_2|\geq2c_K .
  \end{gathered}
  \end{equation}
\end{enumerate}

\noindent\emph{Step 3: a basis adapted to a possible double root.} By point 3, the left pair may collide, and this happens inside
the region of interest: it occurs at the real positive frequency $s$ such that
${|\eta^{1/3}s|^{-3/4}=\kappa_c}$. At such a point
$\{e^{\lambda_1x},e^{\lambda_2x}\}$ is not a basis. Near it, the conditions $y(0)=1$, $y''(0)=0$ written in this basis give
coefficients of the form $\lambda_2^2/\left(\lambda_2^2-\lambda_1^2\right)$ and $-\lambda_1^2/\left(\lambda_2^2-\lambda_1^2\right)$, which blow up although
the solution itself stays bounded. Since $\lambda_2^2-\lambda_1^2=\left(\lambda_2-\lambda_1\right)\left(\lambda_2+\lambda_1\right)$, the cancellation is
built into the basis by dividing by $\lambda_2-\lambda_1$. The condition $\lambda_1\neq\lambda_2$, which fails, is then replaced by
$\lambda_1+\lambda_2\neq0$, which always holds because both roots lie in the left half-plane. The right pair is simple by point 3,
so the usual basis is kept there. Let $\lambda_\tau:=\left(1-\tau\right)\lambda_1+\tau\lambda_2$,
$\tau\in[0,1]$, and
\[
E(x):=\int_0^1x\,e^{\lambda_\tau x}\,d\tau=
\begin{cases}\dfrac{e^{\lambda_2x}-e^{\lambda_1x}}{\lambda_2-\lambda_1}, & \lambda_1\neq\lambda_2,\\[2mm]
x\,e^{\lambda_1x}, & \lambda_1=\lambda_2.\end{cases}
\]
In both cases $\{e^{\lambda_1x},E,e^{\lambda_3\left(x-1\right)},e^{\lambda_4\left(x-1\right)}\}$ is a basis of the solution space of
$\eta y''''+y'+sy=0$, because $\lambda_3\neq\lambda_4$ by point 3 and $\{\lambda_1,\lambda_2\}\cap\{\lambda_3,\lambda_4\}=\emptyset$
by point 2. We have $E(0)=0$ and $E''(0)=\int_0^12\lambda_\tau\,d\tau=\lambda_1+\lambda_2$. Moreover, by convexity,
$\Re\lambda_\tau\leq-c_KN$ and $|\lambda_\tau|\leq C_KN$. Since $x$ appears linearly in the integrand defining $E$,
Leibniz's rule gives $\frac{d^m}{dx^m}\left(x\,e^{\lambda x}\right)=\left(m\lambda^{m-1}+\lambda^mx\right)e^{\lambda x}$. Evaluating at $x=1$
and using $m|\lambda_\tau|^{m-1}\leq mC_K^mN^m$ (recall $C_K\geq1$), we obtain, for $m\leq3$ and $N\geq1$,
\begin{equation}\label{c:Ebounds}
  \begin{gathered}
    |E^{(m)}(1)|\leq\int_0^1\left(m|\lambda_\tau|^{m-1}+|\lambda_\tau|^m\right)e^{\Re\lambda_\tau}d\tau\leq4C_K^mN^me^{-c_KN},\\
    |\lambda_1^me^{\lambda_1}|\leq C_K^mN^me^{-c_KN}.
  \end{gathered}
\end{equation}
In other words, the left modes, together with their derivatives up to order three, are exponentially small at $x=1$:
the derivatives cost powers of $N$, which are absorbed by the factor $e^{-c_KN}$. These bounds are used in Steps 4 and 6.
\\
 
\noindent\emph{Step 4: the layer at $x=1$.} The left modes reach $x=1$ damped by $e^{-c_KN}$, and the conditions at $x=1$ force
the amplitudes of the right modes to be of the same order. Write $y=\alpha e^{\lambda_1x}+\beta E(x)+B_3e^{\lambda_3\left(x-1\right)}+B_4e^{\lambda_4\left(x-1\right)}$ and
$R_m:=\alpha\lambda_1^me^{\lambda_1}+\beta E^{(m)}(1)$. The conditions $y''(1)=y'''(1)=0$ read
$\lambda_3^2B_3+\lambda_4^2B_4=-R_2$ and $\lambda_3^3B_3+\lambda_4^3B_4=-R_3$, hence
\[
B_3=\frac{R_3-\lambda_4R_2}{\lambda_3^2\left(\lambda_4-\lambda_3\right)},\qquad B_4=\frac{R_3-\lambda_3R_2}{\lambda_4^2\left(\lambda_3-\lambda_4\right)} .
\]
By \eqref{c:unifxi} and \eqref{c:Ebounds}, $|R_m|\leq4C_K^mN^m\left(|\alpha|+|\beta|\right)e^{-c_KN}$. Therefore
\begin{equation}\label{c:Bbound}
|B_3|+|B_4|\leq C_B\,\left(|\alpha|+|\beta|\right)\,e^{-c_KN},\qquad C_B:=\frac{16C_K^3}{c_K^2d_K}.
\end{equation}
 
\noindent\emph{Step 5: the layer at $x=0$.} The conditions $y(0)=y_0$ and $y''(0)=0$ read
\[
\alpha+B_3e^{-\lambda_3}+B_4e^{-\lambda_4}=y_0,\qquad
\alpha\lambda_1^2+\beta\left(\lambda_1+\lambda_2\right)+\lambda_3^2B_3e^{-\lambda_3}+\lambda_4^2B_4e^{-\lambda_4}=0 .
\]
Put $\beta=N\tilde\beta$ and divide the second equation by $N^2$. This normalization is natural because
$E''(0)=\lambda_1+\lambda_2$ is of order $N$, while $\lambda_1^2$ is of order $N^2$, so that all the entries of the resulting matrix
are of order one. By Step~4, $B_3,B_4$ are linear in $(\alpha,\beta)$, and
$|e^{-\lambda_j}|\leq e^{-c_KN}$ for $j=3, 4$. We obtain a linear system
\[
\left(A_0+E_N\right)\begin{pmatrix}\alpha\\ \tilde\beta\end{pmatrix}=\begin{pmatrix}y_0\\0\end{pmatrix},\qquad
A_0:=\begin{pmatrix}1&0\\ \xi_1^2&\xi_1+\xi_2\end{pmatrix},
\]
where each entry of the $2\times2$ matrix $E_N$ is bounded by $e_N:=\left(1+C_K^2\right)C_BNe^{-2c_KN}$: each perturbation term carries
one factor $e^{-c_KN}$ from \eqref{c:Bbound} and another one from $e^{-\lambda_j}$. Since
$|\det A_0|=|\xi_1+\xi_2|\geq2c_K$, the entries of $A_0^{-1}$ are bounded by $a_K:=1+\left(C_K^2+1\right)/(2c_K)$. Choose $N_0\geq1$
such that $4a_Ke_N\leq\frac12$ for $N\geq N_0$. In the operator norm induced by the maximum norm,
$\|A_0^{-1}E_N\|\leq4a_Ke_N\leq\frac12$, so $A_0+E_N$ is invertible and $|\alpha|+|\tilde\beta|\leq8a_K|y_0|$. If $y_0=0$
this gives $\alpha=\tilde\beta=0$ and then $B_3=B_4=0$: the homogeneous problem has only the trivial solution. Since $A_0+E_N$ is invertible, the problem is
also solvable for every $y_0$, so $s$ is regular for $\eta$.
\\
 
\noindent\emph{Step 6: the output.} For $y_0=1$, by \eqref{c:Ebounds} with $m=0$: $|E(1)|\leq e^{-c_KN}$ and \eqref{c:Bbound},

\begin{align*}
  |y(1)|&=|\alpha e^{\lambda_1}+N\tilde\beta E(1)+B_3+B_4|\\
  &\leq\left(1+C_B\right)\left(|\alpha|+N|\tilde\beta|\right)e^{-c_KN}\\
  &\leq8a_K\left(1+C_B\right)\,N\,e^{-c_KN}\leq Ce^{-cN},
\end{align*}

with $c:=c_K/2$ and $C:=8a_K\left(1+C_B\right)\sup_{N\geq1}Ne^{-c_KN/2}$.
\end{proof}
 
\begin{remark}\label{c:rem-B}
In physical terms, the signal enters at $x=0$ and is transported by modes decaying like $e^{-c_KNx}$ with $N$ very large,
so that almost nothing reaches $x=1$, and the boundary layers only adjust to this small quantity. The difficulty of the
proof is not the decay itself but the uniformity. No boundary condition may degenerate anywhere in the range of
frequencies, and the only place where this could happen, the collision of the left pair, is handled by Step 3. When $\kappa\to0$, \eqref{c:Qxi} tends to $\xi^4=-e^{i\theta}$. The non-degeneracy that makes the amplitudes of the
growing and decaying modes negligible is exactly \eqref{c:unifxi}, the left pair satisfies $\xi_1+\xi_2\neq0$, and the
right pair is simple and bounded away from $0$.
\end{remark}

We have the following two direct consequences, proved in Appendices \ref{app:cor-B1} and \ref{app:cor-polefree}, respectively.
 
\begin{corollary}[Large frequencies in the regime of Theorem~\ref{thm: stabilitysharpness}]\label{c:cor-B1}
Let $M\geq0$ and $\zeta_*\in(0,1]$. There exist $\eta_B=\eta_B(M,\zeta_*)\in(0,1]$ and $c,C>0$ such that, for every
$\eta\in(0,\eta_B]$ and every $s\in\C$ with $\Re s\geq-M$ and $|\eta^{1/3}s|\geq\zeta_*$, $s$ is regular for $\eta$ and
\[
|f_\eta(s)|\leq C\exp\left(-c\,\eta^{-1/4}|s|^{1/4}\right)\leq C\exp\left(-c\,\zeta_*^{1/4}\eta^{-1/3}\right).
\]
\end{corollary}
 
\begin{corollary}[Pole-free half-planes]\label{c:cor-polefree}
For every $M\geq0$ there is $\eta_P(M)>0$ such that, for every $\eta\in(0,\eta_P(M)]$, every $s$ with $\Re s\geq-M$ is regular for $\eta$. In particular $f_\eta(s)$ and
$f_\eta(s)^2-f_\eta(s)e^{-s}-1$ are defined for all such $s$.
\end{corollary}
 
\subsection{The upper bound}\label{c:sec-upper}
 Before studying sequences of roots of \eqref{eq: fourthorderCharEq}, we state two auxiliary results which will be used in the proof.

\begin{lemma}\label{c:lem-lowerF}
Let $M\geq0$, $s\in\C$ with $\Re s\geq-M$, and $F\in\C$ with $F^2-Fe^{-s}-1=0$. Then
$|F|\geq m_M:=\frac12\left(\sqrt{e^{2M}+4}-e^M\right)>0$.
\end{lemma}
 
\begin{proof}
See Appendix~\ref{app:lem-lowerF}.
\end{proof}
 
\begin{lemma}[Damping of the finite branch]\label{c:lem-damping}
Decreasing $\zeta_0$ if necessary, $\Re\lambda_1(s)\leq-\frac12\Re s$ for every $\eta>0$ and every $s\in\C$ with $\Re s\geq0$
and $|\eta^{1/3}s|\leq\zeta_0$.
\end{lemma}
 
\begin{proof}
See Appendix~\ref{app:lem-damping}.
\end{proof}
 
Throughout this subsection, we consider a sequence 
\begin{equation}\label{eq: visc}
  (\eta_n)_{n\in\N}\subset \R, \quad \eta_n>0, \quad \eta_n\to 0^+,
\end{equation}
and an associated sequence 
\begin{equation}\label{eq: poles}
  (s_n)_{n\in \N}\subset \C, \quad s_n = \sigma_n + i \omega_n, \quad s_n\in \Sp_{\eta_n}, \ \forall n\in \N.
\end{equation}
We use accumulation
points of the real sequence $(\sigma_n)_{n\in\N}$ in the extended line $\left[-\infty,+\infty\right]$.
 
\begin{lemma}\label{c:lem-32}
Let $(\eta_n)_{n\in\N}$ be a sequence of the form \eqref{eq: visc} and $(s_n)_{n\in\N}$ be any associated sequence of the form \eqref{eq: poles}. Then:
\begin{enumerate}[itemsep=1pt]
\item $|s_n|\to+\infty$. In other words, $(s_n)_{n\in\N}$ has no bounded subsequence.
\item $(\sigma_n)_{n\in\N}$ is bounded from above.
\item If $\sigma_n\to\bar\sigma\in\R$ along a subsequence, then along it $\zeta_n\to0$ and $|\omega_n|\to+\infty$.
\end{enumerate}
\end{lemma}
 
\begin{proof}
  \begin{enumerate}
    \item Suppose that a subsequence converges to some $\bar s\in\C$, and let $K$ be the closed unit disc centered at $\bar s$.
By Corollary~\ref{c:cor-compact}, $\sup_K|f_{\eta_n}-e^{-s}|\to0$, hence $f_{\eta_n}(s_n)\to e^{-\bar s}$. Passing to the
limit in \eqref{eq: fourthorderCharEq} gives $e^{-2\bar s}-e^{-2\bar s}-1=0$, a contradiction.

    \item For $\Re s$ large, both $e^{-s}$ and the transmitted signal $f_\eta(s)\approx e^{\lambda_1(s)}$ are small, so that
\[f_\eta(s)\left(f_\eta(s)-e^{-s}\right)=1\] 
cannot hold. For the sake of contradiction, assume that $\sigma_n\to+\infty$ along a subsequence. Along it
$\sigma_n\geq0$, so $|f_{\eta_n}(s_n)|\geq m_0$ by Lemma~\ref{c:lem-lowerF}. By Corollary~\ref{c:cor-B1} with $M=0$ and
$\zeta_*=\zeta_0$, $|\zeta_n|<\zeta_0$ for all large $n$, since otherwise
$|f_{\eta_n}(s_n)|\leq C\exp\left(-c\zeta_0^{1/4}\eta_n^{-1/3}\right)\to0$ along a further subsequence. Hence, for large $n$,
Proposition~\ref{c:prop-A} and Lemma~\ref{c:lem-damping} apply and, since $|\Phi|\leq2$ and $|r_2|\leq\frac12$,
\[
|f_{\eta_n}(s_n)|\leq3\,e^{\Re\lambda_1(s_n)}+C_Ae^{-\eta_n^{-1/3}/2}\leq3\,e^{-\sigma_n/2}+C_Ae^{-\eta_n^{-1/3}/2}\longrightarrow0,
\]
a contradiction.

    \item Let $\sigma_n\to\bar\sigma\in\R$ along a subsequence. Along it $\sigma_n\geq-M:=\bar\sigma-1$ for large $n$, so
$|f_{\eta_n}(s_n)|\geq m_M$ by Lemma~\ref{c:lem-lowerF}. For every $\zeta_*\in(0,\zeta_0]$, Corollary~\ref{c:cor-B1} gives
$|\zeta_n|<\zeta_*$ for large $n$ (otherwise $f_{\eta_n}(s_n)\to0$ along a further subsequence). Hence $\zeta_n\to0$.
Finally $|\omega_n|\to\infty$ by item~1, since $\sigma_n$ is bounded.
\end{enumerate}
 
\end{proof}
 
\begin{lemma}\label{c:lem-33}
Let $(\eta_n)_{n\in\N}$ be a sequence of the form \eqref{eq: visc} and $(s_n)_{n\in\N}$ be any associated sequence of the form \eqref{eq: poles}.
Write $\zeta_n:=\eta_n^{1/3}s_n$ and $t_n:=\eta_n\omega_n^4$.Assume that
$\omega_n\to+\infty$ and $t_n\to\bar t\in\left[0,+\infty\right]$. If $\bar\sigma\in(-\infty,+\infty]$ is an accumulation point of
$(\sigma_n)_{n\in\N}$, then $\bar\sigma\in\R$, $\bar t\in\left(0,+\infty\right)$ and $\bar\sigma=g(\bar t)$, where
\[
g(t):=-\frac t2+\frac12\ln\left(1-e^{-t}\right),\qquad t>0 .
\]
In particular: if $\bar t\in\{0,+\infty\}$, then $\sigma_n\to-\infty$. If $\bar t\in\left(0,+\infty\right)$, then every accumulation
point of $(\sigma_n)_{n\in\N}$ belongs to $\{-\infty,\,g(\bar t)\}$.
\end{lemma}
 
\begin{proof}
Pass to a subsequence with $\sigma_n\to\bar\sigma$. By item 2 and 3 of Lemma~\ref{c:lem-32}, we have $\bar\sigma<+\infty$ and
$\zeta_n\to0$ respectively. Hence Proposition~\ref{c:prop-A} applies for large $n$:
\[
f_{\eta_n}(s_n)=e^{\lambda_{1,n}}\Phi(\zeta_n)\left(1+r_{2,n}\right)+r_{3,n},\qquad \lambda_{1,n}=-s_n-\eta_ns_n^4\left(1+r_{1,n}\right),
\]
with $\Phi(\zeta_n)\to1$, $r_{2,n}, r_{3,n}\to0$ and $|r_{1,n}|\leq C_0|\zeta_n|^3\to0$. Since
$\left(\sigma+i\omega\right)^4=\omega^4\left(1-i\sigma/\omega\right)^4$ and $\sigma_n/\omega_n\to0$,
\begin{equation}\label{c:es4}
\eta_ns_n^4\left(1+r_{1,n}\right)=t_n\left(1+o(1)\right).
\end{equation}
 
\noindent\emph{Case $\bar t<\infty$.} By \eqref{c:es4}, $\lambda_{1,n}+s_n\to-\bar t$, and since $e^{\sigma_n}$ is bounded,
$e^{s_n}f_{\eta_n}(s_n)=e^{-\bar t}+o(1)$. Multiplying \eqref{eq: fourthorderCharEq} by $e^{2s_n}$ gives $\left(e^{s_n}f_{\eta_n}(s_n)\right)^2-e^{s_n}f_{\eta_n}(s_n)=e^{2s_n}$. Taking
moduli and passing to the limit,
\[
e^{2\bar\sigma}=\big|e^{-2\bar t}-e^{-\bar t}\big|=e^{-\bar t}\left(1-e^{-\bar t}\right).
\]
If $\bar t=0$ the right-hand side vanishes, which is impossible. If $\bar t\in\left(0,\infty\right)$, then
$\bar\sigma=\frac12\ln\left(e^{-\bar t}\left(1-e^{-\bar t}\right)\right)=g(\bar t)$.
\\
 
\noindent\emph{Case $\bar t=+\infty$.} By \eqref{c:es4}, $\Re\left(\eta_ns_n^4\left(1+r_{1,n}\right)\right)\geq t_n/2$ for large $n$, so
$\Re\lambda_{1,n}\leq-\sigma_n-t_n/2\to-\infty$, because $\sigma_n$ is bounded below. Therefore
$|f_{\eta_n}(s_n)|\leq3e^{\Re\lambda_{1,n}}+|r_{3,n}|\to0$, which contradicts Lemma~\ref{c:lem-lowerF} with $M$ such that $\sigma_n\geq-M$.
\end{proof}
 
\begin{lemma}\label{c:lem-34}
Every accumulation point $\bar\sigma\in\left[-\infty,+\infty\right]$ of $(\sigma_n)_{n\in\N}$ satisfies $\bar\sigma\leq-\ln(2)$.
\end{lemma}
 
\begin{proof}
If $\bar\sigma=-\infty$ there is nothing to prove, and $\bar\sigma=+\infty$ is excluded by the second point of Lemma~\ref{c:lem-32}. Let
$\bar\sigma\in\R$ and $\sigma_n\to\bar\sigma$ along a subsequence. By Lemma~\ref{c:lem-32}(3), $|\omega_n|\to\infty$. Since
$\Sp_{\eta_n}$ is invariant under conjugation, see Remark~\ref{c:rem-spectrum}. Replacing $s_n$ by $\bar s_n$ when $\omega_n<0$,
which changes neither $\sigma_n$ nor $t_n$, we may assume $\omega_n\to+\infty$. Extract $t_n\to\bar t\in\left[0,\infty\right]$, by item 3 of 
Lemma~\ref{c:lem-33}, $\bar t\in\left(0,\infty\right)$ and $\bar\sigma=g(\bar t)$.
\\
 
\noindent It remains to bound $g$ on $\left(0,\infty\right)$. We have
\[
g'(t)=-\frac12+\frac{1}{2\left(e^t-1\right)},\qquad g''(t)=-\frac{e^t}{2\left(e^t-1\right)^2}<0 .
\]
Thus $g$ is concave and $g'$ is strictly decreasing, with $g'(t)\to+\infty$ as $t\to0^+$, $g'(t)\to-\frac12$ as $t\to\infty$,
and a unique zero at $t^*=\ln(2)$ where $e^{t^*}-1=1$. Hence $g$ is increasing on $(0,\ln(2)]$ and decreasing on
$[\ln(2),\infty)$, and $g(t)\leq g\left(\ln(2)\right)=-\frac{\ln(2)}2+\frac12\ln\frac12=-\ln(2)$ for every $t>0$.
\end{proof}
 
\begin{proof}[Proof of the upper bound in Theorem~\ref{thm: stabilitysharpness}]
Suppose that the upper bound fails for some $\delta>0$. Then there are $\eta_n\in\left(0,1/n\right)$ and $s_n\in\Sp_{\eta_n}$ with
$\sigma_n>-\ln(2)+\delta$ for all $n$. The real sequence $(\sigma_n)_{n\in\N}$ has an accumulation point
$\bar\sigma\in\left[-\infty,+\infty\right]$, and necessarily $\bar\sigma\geq-\ln(2)+\delta>-\ln(2)$. This contradicts
Lemma~\ref{c:lem-34}. Hence there is $\tilde\eta_1(\delta)>0$ such that $\sigma_\eta\leq-\ln(2)+\delta$ for all
$\eta\in\left(0,\tilde\eta_1(\delta)\right)$.
\end{proof}
 
\begin{remark}[The threshold $\eta_1$]\label{c:rem-etaup}
As in \cite{BCH2025Viscosity}, we denote by $\eta_1$ the threshold given by the upper bound. Precisely, we set $\eta_1(\delta):=\tilde\eta_1(\delta)$. Thus $\sigma_\eta\leq-\ln(2)+\delta$ for every
$\eta\in\left(0,\eta_1(\delta)\right)$.
\end{remark}
 
\subsection{The lower bound}\label{c:sec-lower}
 
\medskip\noindent\emph{{Step 1: rescaled variables and the critical point.}} Let $\gamma:=\eta^{1/4},\varphi:= \gamma \omega$ and write $s=\sigma+i\varphi/\gamma$ with $\omega=\Im s$. The critical point identified in Lemma~\ref{c:lem-34} is $t^*=\ln(2)$,
$\sigma_*=g\left(t^*\right)=-\ln(2)$. In the new variables, $(\sigma_*,\varphi_*)=\left(-\ln(2),\left(\ln 2\right)^{1/4}\right)$.
 
\medskip\noindent\emph{{Step 2: factorization into a fast phase and a slow amplitude.}} Let $\delta>0$ and fix $\delta_0\in(0,\delta)$.
Let $Q:=\left[\sigma_*-\delta_0,\sigma_*+\delta_0\right]\times\left[\varphi_*/2,2\varphi_*\right]$ and, for $\eta=\gamma^4\leq\eta_Q$
given by Corollary~\ref{c:cor-resonant}, define on $Q$
\[
h_\gamma(\sigma,\varphi):=e^{2s}\left(f_\eta(s)^2-f_\eta(s)e^{-s}\right),\qquad s=\sigma+i\varphi/\gamma,
\qquad h_0(\varphi):=e^{-2\varphi^4}-e^{-\varphi^4}<0 .
\]
The function $h_\gamma$ is well defined and continuous on $Q$ by Corollary~\ref{c:cor-resonant}, every
$s=\sigma+i\varphi/\gamma$ with $(\sigma,\varphi)\in Q$ is regular for $\eta$ and satisfies $|\eta^{1/3}s|\leq\zeta_0$, and $f_\eta$
is continuous on this set by Proposition~\ref{c:prop-A}.
 
\begin{lemma}\label{c:lem-hgamma}
There is $C_2>0$ such that $|h_\gamma(\sigma,\varphi)-h_0(\varphi)|\leq C_2\,\eta^{1/6}$ for all $(\sigma,\varphi)\in Q$ and
$\eta=\gamma^4\in(0,\eta_Q]$.
\end{lemma}
 
\begin{proof}
 See Appendix~\ref{app:hgamma}.
\end{proof}
 
Let $\mathcal G(s):=f_\eta(s)^2-f_\eta(s)e^{-s}-1$ for $s$ regular for $\eta$, so that the zeros of
$\mathcal G$ are exactly the points of $\Sp_\eta$. With $e^{-2s}=e^{-2\sigma}e^{-2i\varphi/\gamma}$,
\[
\mathcal G(s)=e^{-2\sigma}e^{-2i\varphi/\gamma}h_\gamma(\sigma,\varphi)-1
=-\left(e^{-2i\varphi/\gamma}\cA_\gamma(\sigma,\varphi)+\cB\right),\qquad
\cA_\gamma:=-e^{-2\sigma}h_\gamma,\quad \cB:=1,
\]
and $\cA_0(\sigma,\varphi):=-e^{-2\sigma}h_0(\varphi)=e^{-2\sigma}\left(e^{-\varphi^4}-e^{-2\varphi^4}\right)>0$. With the canonical
isomorphism $\Psi\left(x+iy\right):=(x,y)$ we set $\cH_\gamma:=\Psi\left(e^{-2i\varphi/\gamma}\cA_\gamma+\cB\right)$. Note that the zeros of
$\cH_\gamma$ are exactly the zeros of $\mathcal G$.
 
\medskip\noindent\emph{{Step 3: the structural hypotheses (H1)--(H2) from \cite[Remark 4]{BCH2025Viscosity}}.} With $t_*=\varphi_*^4=\ln(2)$,\\

\noindent\emph{(H1) Amplitude matching at the critical point.}
\[
\cA_0(\sigma_*,\varphi_*)=e^{2\ln 2}\left(\frac12-\frac14\right)=1=\cB.\]    
\noindent\emph{(H2) Strict monotonicity of the modulus gap.} The map
\[\sigma\mapsto|\cA_0(\sigma,\varphi_*)|-|\cB|=e^{-2\left(\sigma-\sigma_*\right)}-1,
\] is strictly decreasing, positive for $\sigma<\sigma_*$ and negative for $\sigma>\sigma_*$.
 
\medskip\noindent\emph{{Step 4: the resonant rectangle and boundary estimates.}} Let
\[
k(\gamma):=\left\lfloor\frac{\varphi_*}{\pi\gamma}\right\rfloor,\qquad
\Omega_\gamma:=\left(\sigma_*-\delta_0,\sigma_*+\delta_0\right)\times\left(k(\gamma)\pi\gamma,\left(k(\gamma)+1\right)\pi\gamma\right).
\]
Then $|\varphi-\varphi_*|\leq\pi\gamma$ on $\overline{\Omega}_\gamma$, and $\overline{\Omega}_\gamma\subset Q$ as soon as
$\pi\gamma\leq\varphi_*/2$. Let $a(\sigma):=e^{-2\left(\sigma-\sigma_*\right)}=\cA_0(\sigma,\varphi_*)$ and $\eps_\gamma:=\cA_\gamma/a-1$, so that
$\cA_\gamma=a\left(1+\eps_\gamma\right)$.
\\
 
We claim that $\sup_{\overline\Omega_\gamma}|\eps_\gamma|\leq C_1\eta^{1/6}$, with $C_1:=4e^{4\delta_0}\left(C_2+C_4\right)$ and
$C_4:=\pi\sup_{\left[\varphi_*/2,2\varphi_*\right]}|h_0'|$. Indeed,
$\cA_\gamma-a=-e^{-2\sigma}\left(h_\gamma-h_0(\varphi)\right)-e^{-2\sigma}\left(h_0(\varphi)-h_0(\varphi_*)\right)$, where
$e^{-2\sigma}\leq K_1:=e^{-2\sigma_*+2\delta_0}=4e^{2\delta_0}$ on $\overline\Omega_\gamma$. The first term is bounded by
$K_1C_2\eta^{1/6}$ thanks to Lemma~\ref{c:lem-hgamma}. By the mean value theorem, with
$h_0'(\varphi)=4\varphi^3\left(e^{-\varphi^4}-2e^{-2\varphi^4}\right)$, the second is bounded by $K_1C_4\gamma\leq K_1C_4\eta^{1/6}$, since
$\gamma=\eta^{1/4}\leq\eta^{1/6}$. Dividing by $a\geq e^{-2\delta_0}$ gives the claim.
\\
 
\begin{figure}[h]
\centering
\includegraphics[width=1.05\textwidth]{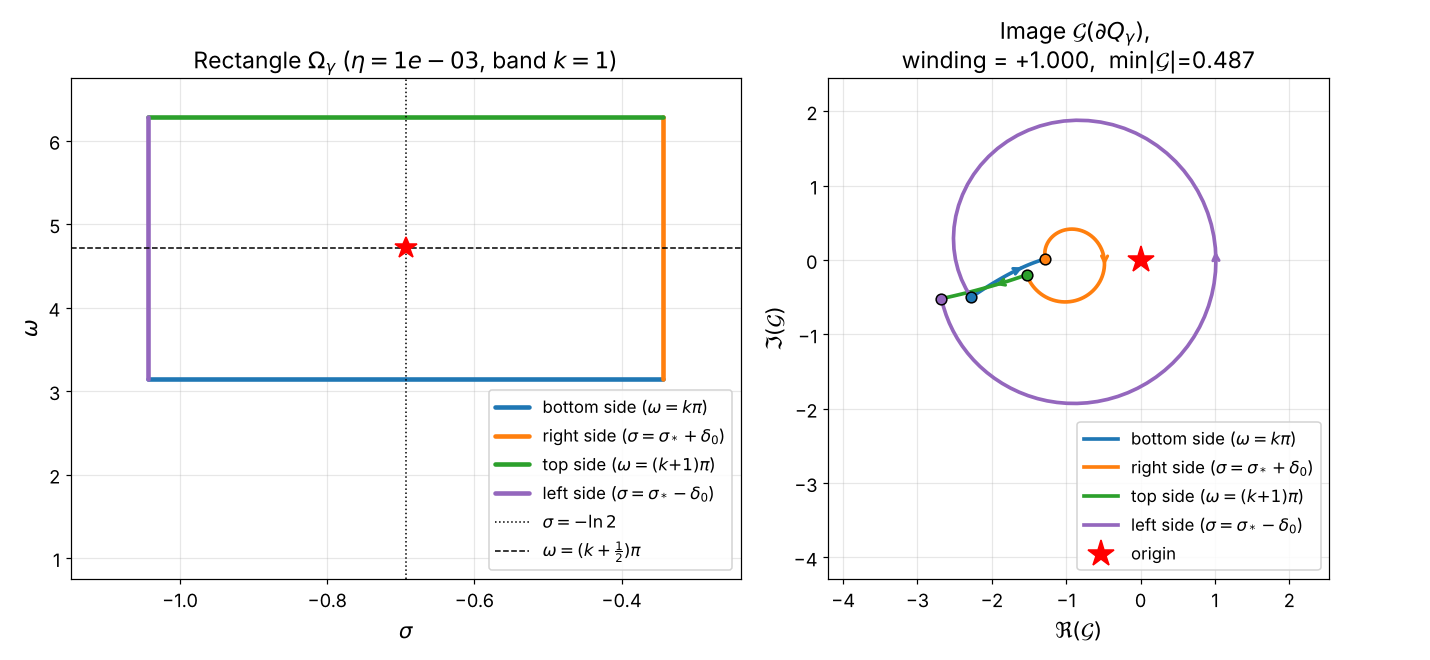}
\caption{Image of the rectangle $\Omega_\gamma$ under $\mathcal{G}$.}
\label{fig: chap5_windingnumber}
\end{figure}

Let $c(\delta_0):=\frac12\min\{1-e^{-2\delta_0},e^{2\delta_0}-1\}=\frac12\left(1-e^{-2\delta_0}\right)$. We fix $\gamma_1>0$ such that, for
every $\gamma\in(0,\gamma_1]$,
\begin{equation}\label{c:gamma1}
\gamma^4\leq\eta_Q,\qquad \pi\gamma\leq\varphi_*/2,\qquad C_1\gamma^{2/3}\leq\tfrac12c(\delta_0)e^{-2\delta_0}.
\end{equation}
For $\mu\in[0,1]$ put $\cA^{(\mu)}:=\left(1-\mu\right)a+\mu\cA_\gamma=a\left(1+\mu\eps_\gamma\right)$, and note that
$|\mu\eps_\gamma|\leq|\eps_\gamma|\leq\frac12c(\delta_0)e^{-2\delta_0}<\frac12$. Then, for $\gamma\in(0,\gamma_1]$ and every
$\mu\in[0,1]$:
\begin{enumerate}[label=\textup{(\alph*)},itemsep=1pt]
\item \emph{Horizontal sides}, where $\varphi=k(\gamma)\pi\gamma$ or $\varphi=\left(k(\gamma)+1\right)\pi\gamma$. There
$e^{-2i\varphi/\gamma}=1$, so the first component of $\Psi\left(e^{-2i\varphi/\gamma}\cA^{(\mu)}+1\right)$ is
$\Re\cA^{(\mu)}+1\geq a\left(1-|\eps_\gamma|\right)+1>1$.
\item \emph{Vertical sides.} On $\sigma=\sigma_*+\delta_0$: $a=e^{-2\delta_0}=1-2c(\delta_0)$ and
$|\cA^{(\mu)}|\leq a+a|\eps_\gamma|\leq1-2c(\delta_0)+\frac12c(\delta_0)<1-c(\delta_0)$. On $\sigma=\sigma_*-\delta_0$:
$a=e^{2\delta_0}\geq1+2c(\delta_0)$ and $|\cA^{(\mu)}|\geq a-a|\eps_\gamma|\geq e^{2\delta_0}-\frac12c(\delta_0)>1+c(\delta_0)$.
In both cases $|e^{-2i\varphi/\gamma}\cA^{(\mu)}+1|\geq\big||\cA^{(\mu)}|-1\big|>c(\delta_0)$.
\end{enumerate}
Consequently $\Psi\left(e^{-2i\varphi/\gamma}\cA^{(\mu)}+1\right)\neq0$ on $\partial\Omega_\gamma$ for every $\mu\in[0,1]$. For $\mu=1$
this is $\cH_\gamma\neq0$ on $\partial\Omega_\gamma$. Figure \ref{fig: chap5_windingnumber} illustrates this step for $\eta=10^{-3}$. It shows the rectangle $\Omega_\gamma$ and its image under $\mathcal{G}$, which winds once around the origin.
 
\medskip\noindent\emph{{Step 5: the model map and its degree.}} Let $\cM_\gamma(\sigma,\varphi):=\Psi\left(e^{-2i\varphi/\gamma}a(\sigma)+1\right)$, the case $\mu=0$. A point of $\overline\Omega_\gamma$ is a zero of $\cM_\gamma$ if and only if $e^{-2i\varphi/\gamma}a(\sigma)=-1$.
Since $a(\sigma)>0$ and $|e^{-2i\varphi/\gamma}|=1$, this means $a(\sigma)=1$, i.e.\ $\sigma=\sigma_*$, and
$e^{-2i\varphi/\gamma}=-1$, i.e.\ $\varphi/\gamma\in\left(\Z+\frac12\right)\pi$. As $\varphi/\gamma\in\left[k(\gamma)\pi,\left(k(\gamma)+1\right)\pi\right]$ on
$\overline\Omega_\gamma$, the only zero is
\[
P_*:=\left(\sigma_*,\left(k(\gamma)+\tfrac12\right)\pi\gamma\right)\in\Omega_\gamma .
\]
Writing $\cM_\gamma=(M_1,M_2)$ with $M_1=a(\sigma)\cos\left(2\varphi/\gamma\right)+1$, $M_2=-a(\sigma)\sin\left(2\varphi/\gamma\right)$, and using
$a'=-2a$, $a(\sigma_*)=1$ and $2\varphi/\gamma=\left(2k(\gamma)+1\right)\pi$ at $P_*$, one finds
\[
J(\cM_\gamma)(P_*)=\begin{pmatrix}2&0\\0&2/\gamma\end{pmatrix},\qquad \det J(\cM_\gamma)(P_*)=\frac4\gamma>0 .
\]
Since $\cM_\gamma\neq0$ on $\partial\Omega_\gamma$ by Step~4 and $P_*$ is its only zero, which is regular, the
regular-value formula \cite[Theorem~B.1, (B.4)]{Coron2007} gives
\[
\deg(\cM_\gamma,\Omega_\gamma,0)=\sgn\det J(\cM_\gamma)(P_*)=+1 .
\]
 
\medskip\noindent\emph{{Step 6: homotopy and transfer of the degree.}} Consider
\[
\cH_\gamma^{(\mu)}(\sigma,\varphi):=\Psi\left(e^{-2i\varphi/\gamma}\left[\left(1-\mu\right)a(\sigma)+\mu\cA_\gamma(\sigma,\varphi)\right]+1\right)
=\Psi\left(e^{-2i\varphi/\gamma}\cA^{(\mu)}+1\right),\qquad\mu\in[0,1],
\]
so that $\cH^{(0)}_\gamma=\cM_\gamma$ and $\cH^{(1)}_\gamma=\cH_\gamma$. The map $(\mu,\sigma,\varphi)\mapsto\cH_\gamma^{(\mu)}(\sigma,\varphi)$
is continuous on $[0,1]\times\overline\Omega_\gamma$, because $f_\eta$ is continuous at every point $\sigma+i\varphi/\gamma$ with
$(\sigma,\varphi)\in\overline\Omega_\gamma\subset Q$ (Step~2). By Step~4,
$\cH^{(\mu)}_\gamma\neq0$ on $\partial\Omega_\gamma$ for all $\mu\in[0,1]$. By the homotopy invariance of the degree
\cite[Proposition~B.8]{Coron2007},
\[
\deg(\cH_\gamma,\Omega_\gamma,0)=\deg(\cM_\gamma,\Omega_\gamma,0)=+1\neq0,
\]
and by \cite[Proposition~B.10]{Coron2007}, $\cH_\gamma$ has a zero $(\tilde\sigma,\tilde\varphi)\in\Omega_\gamma$.
 
\medskip\noindent\emph{{Step 7: conclusion.}} Let $\gamma\in(0,\gamma_1]$, $\eta=\gamma^4$ and $\tilde s:=\tilde\sigma+i\tilde\varphi/\gamma$. Then
$\tilde s$ is regular for $\eta$ and $\mathcal G(\tilde s)=0$, so $\tilde s\in\Sp_\eta$, and
$\Re\tilde s=\tilde\sigma>-\ln(2)-\delta_0>-\ln(2)-\delta$. Hence $\sigma_\eta>-\ln(2)-\delta$ for every $\eta\in(0,\gamma_1^4]$.
Moreover $\eta(\Im\tilde s)^4=\tilde\varphi^4$ and $|\tilde\varphi-\varphi_*|<\pi\gamma$, so $\eta(\Im\tilde s)^4\to\ln(2)$ as
$\eta\to0^+$. Together with Section~\ref{c:sec-upper}, \eqref{eq: fourthordersharpness2} holds with
\[
\eta_*:=\min\{\eta_1(\delta),\gamma_1^4\},
\]
which completes the proof of Theorem~\ref{thm: stabilitysharpness}. \qed

\section{Proof of Theorem \ref{thm:robustness}}\label{sec: robustness}
 
Recall the perturbed system \eqref{eq: fourthordersystepsilon1}, its single-channel transfer function
$f_{\eta,\varepsilon}(s)$, defined at every point $s$ that is regular for $(\eta,\eps)$
by Definition~\ref{c:def-regular}, and the closed-loop characteristic equation \eqref{eq:char-eq-eps1}. The spectrum
$\mathcal S_{\eta,\varepsilon}$ and the spectral abscissa $\sigma_{\eta,\varepsilon}$ are those of
Definition~\ref{def: spectrumeps}, with the convention of Remark~\ref{c:rem-spectrum}. As before, the loop is closed with the
dead-beat feedback $U(t)=-Y\left(t-1\right)$, which uses the nominal delay $\tau=1$ and thus ignores the uncertainty $\varepsilon$. We
assume that the same perturbation $\varepsilon$ affects both channels.
\\
 
The perturbation of the velocity can be absorbed by a change of the parameters. Indeed, let $\varepsilon>-1$. Dividing the equation in \eqref{c:bvp} by $1+\varepsilon$ gives $\frac{s}{1+\varepsilon}\,y+y'+\frac{\eta}{1+\varepsilon}\,y''''=0$, with the same boundary conditions, so for every $y_0$ both problems have the same solutions. Hence $s$ is regular for $(\eta,\varepsilon)$ if and only if $s/(1+\varepsilon)$ is regular for $\eta/(1+\varepsilon)$, and in that case,
\begin{equation}\label{c:scaling}
f_{\eta,\eps}(s)=f_{\eta/(1+\eps)}\left(\frac{s}{1+\eps}\right).
\end{equation}

By \eqref{c:scaling}, the perturbation creates no new boundary layer: it only changes the parameter $\eta$ and rescales the
frequency, and the results of Section~\ref{sec: asymptotic} apply with $\eps=0$.
 
\begin{proof}[Proof of Theorem~\ref{thm:robustness}]
Let $\delta>0$, put  $\eta^*(\delta):=\min\{\eta_1(\delta),\eta_P(1)\}$, with $\eta_1(\delta)$ from Remark~\ref{c:rem-etaup} and
$\eta_P(1)$ from Corollary~\ref{c:cor-polefree}, and let $\eta\in\left(0,\eta^*(\delta)\right)$.
Then $\sigma_{\eta,0}\leq-\ln(2)+\delta$, and every $s$ with $\Re s\geq-1$ is regular for $\eta$. Suppose,
towards a contradiction, that no $\eps_1>0$ satisfies \eqref{eq: fourthorderrobustness}. Then there are $\eps_n\to0$ of either sign, with
$|\eps_n|\leq\frac12$ and $s_n\in\Sp_{\eta,\eps_n}$ with
\begin{equation}\label{c:47}
\Re s_n>-\ln(2)+2\delta\qquad\text{for all }n .
\end{equation}
By Lemma~\ref{c:lem-lowerF} with $M=\ln 2$, $|f_{\eta,\eps_n}(s_n)|\geq m_{\ln 2}>0$ for all $n$.
\\
 
\noindent\emph{{Case 1: a bounded subsequence.}} Extract $s_n\to\bar s\in\C$. By \eqref{c:47},
$\Re\bar s\geq-\ln(2)+2\delta>-1$, so $\bar s$ is regular for $\eta$. By Proposition \ref{c:prop-A},
$f_{\eta,\eps_n}(s_n)\to f_\eta(\bar s)$, and passing to the limit in \eqref{eq:char-eq-eps1} gives
$f_\eta(\bar s)^2-f_\eta(\bar s)e^{-\bar s}-1=0$, i.e.\ $\bar s\in\Sp_{\eta,0}=\Sp_\eta$. Hence
$\Re\bar s\leq\sigma_{\eta,0}\leq-\ln(2)+\delta<-\ln(2)+2\delta\leq\Re\bar s$, a contradiction.
\\
 
\noindent\emph{{Case 2: $|s_n|\to\infty$.}} Write $s_n=r_ne^{i\theta_n}$ with $\theta_n\in(-\pi,\pi]$. By \eqref{c:47},
$\cos\theta_n>\left(-\ln(2)+2\delta\right)/r_n\geq-\ln(2)/r_n\to0$. Hence, with $\theta_1$ given by Proposition~\ref{c:prop-B} for
$\kappa_0:=1$, we have $|\theta_n|\leq\frac\pi2+\theta_1$ for $n$ large. Put $\eta_n:=\eta/(1+\eps_n)$ and
$\tilde s_n:=s_n/(1+\eps_n)$, so that $\arg\tilde s_n=\theta_n$ and the quantities \eqref{c:Nkappa} associated with
$(\eta_n,\tilde s_n)$ are
\[
\tilde N_n:=\left(\frac{|\tilde s_n|}{\eta_n}\right)^{1/4}=\left(\frac{r_n}{\eta}\right)^{1/4}\longrightarrow\infty,\qquad
\tilde\kappa_n:=|\eta_n^{1/3}\tilde s_n|^{-3/4}=\frac{1+\eps_n}{\eta^{1/4}r_n^{3/4}}\longrightarrow0 .
\]
Thus, for $n$ large, we apply Proposition \ref{c:prop-A} before using Proposition \ref{c:prop-B} for $(\eta_n,\tilde s_n)$ with $\kappa_0=1$, obtaining that $\tilde s_n$ is regular for
$\eta_n$ and
\[
|f_{\eta,\eps_n}(s_n)|=|f_{\eta_n}(\tilde s_n)|\leq Ce^{-c\tilde N_n}\longrightarrow0 .
\]

\noindent This contradicts $|f_{\eta,\eps_n}(s_n)|\geq m_{\ln 2}$. 

 
\noindent Both cases lead to a contradiction, which proves the theorem.
\end{proof}

\section{Conclusion}\label{sec: conclusion}

We have studied the output-feedback stabilization of an unstable open-loop system consisting of two interconnected
transport equations perturbed by fourth-order terms, with anti-located control and observation. We have shown that,
even without an explicit description of the roots of the quartic characteristic equation, stability can be analyzed
through an asymptotic study of these roots and of the transfer function.
\\

In the presence of a small fourth-order term, the system is robustly stabilized by a delayed output feedback. The optimality of the margin $-\ln(2)$ follows from a topological degree argument at the
resonant frequencies $\omega\sim\eta^{-1/4}$ created by the fourth-order term, and this margin is robust to small
perturbations of the transport velocity.
\\

The significance of this extension for control is not merely that the same numerical value $-\ln(2)$ reappears, but that
it does so through a mechanism that does not rely on explicit root formulas. This indicates that the mechanism by which
diffusion restores robustness is not an artifact of the specific second-order model first studied in
\cite{BCH2022,BCH2025Viscosity}, but persists under a structurally different, higher-order dissipative perturbation.
\\

Two aspects of the present analysis should be stated as explicit limitations. First, as noted in
Remark~\ref{rem: SMT}, deducing exponential decay of the full closed-loop system (and not only focusing on the output) from $\sigma_\eta<0$ requires a spectral
mapping (or spectrum-determined growth) property. Second, we have only treated the specific dead-beat gain
$U(t)=-Y\left(t-1\right)$ from \cite{BCH2022}. Whether the sharp margin $-\ln(2)$ and its robustness persist for general feedback
gains remains open.
\\

Natural directions for future work include determining whether the mechanism identified here persists for higher-order
dissipative regularizations $\left(-1\right)^m\eta\,\partial_x^{2m}$, $m\geq3$ with the sign convention of
\eqref{eq: fourthordersyst} and the appropriate choice of boundary conditions. {Our computations indicate that the strategy developed here seems to work in the general case, at the cost of heavier algebra.} Along the same line, another possible extension is the regularization of the second
equation of the open-loop system of \cite{BCH2025Viscosity} with a fourth-order dissipative term. If this is the case, the natural extension is to consider a general diffusion term and study the stabilization properties for the system.


\section*{Appendix}
\appendix

\section{Proof of Lemma \ref{c:lem-wp}}\label{app:lem-wp}

1. The solutions of the differential equation in \eqref{c:bvp} form a vector space $V$ of dimension four, since it is a
linear equation of order four with constant coefficients.  Consider the linear map $L:V\to\C^4$,
$Ly:=\left(y(0),y''(0),y''(1),y'''(1)\right)$. If $s$ is regular, uniqueness for $y_0=0$ is part of Definition~\ref{c:def-regular}. Conversely, if \eqref{c:bvp} with $y_0=0$ has only the trivial solution, then $L$ is injective. An injective linear map between
spaces of the same finite dimension is bijective, so for every $y_0$ there is exactly one $y\in V$ with
$Ly=\left(y_0,0,0,0\right)$, that is, $s$ is regular.
\\

\noindent 2. Since $\eta$ and $\eps$ are real, $y$ solves \eqref{c:bvp} for $s$ and $y_0$ if and only if $\bar y$ solves it for $\bar s$
and $\bar y_0$. Hence $y\mapsto\bar y$ is a bijection between the solutions at $s$ with datum $y_0$ and the solutions at $\bar s$ with datum $\bar y_0$, so existence and uniqueness for every datum hold at $s$ if and only if they hold at $\bar s$,
and taking $y_0=1$ gives $f_{\eta,\eps}(\bar s)=\overline{f_{\eta,\eps}(s)}$.
\\

\noindent 3. Let $y$ solve \eqref{c:bvp} with $s=0$ and $y_0=0$. Multiply the equation by $\bar y$, integrate over $(0,1)$ and
integrate by parts twice in the last term:
\[
\int_0^1y''''\,\bar y\,dx=\left[y'''\,\bar y-y''\,\bar y'\right]_0^1+\int_0^1|y''|^2dx=\int_0^1|y''|^2dx,
\]
because $y(0)=y''(0)=y''(1)=y'''(1)=0$. Since $2\Re\int_0^1y'\,\bar y\,dx=|y(1)|^2$, the real part of the resulting identity is
$\frac{1+\eps}2|y(1)|^2+\eta\int_0^1|y''|^2dx=0$. Hence $y''\equiv0$ and $y(1)=0$, which together with $y(0)=0$ gives $y\equiv0$.
So $s=0$ is regular by point 1, and since $y\equiv1$ solves \eqref{c:bvp} with $s=0$ and $y_0=1$, $f_{\eta,\eps}(0)=1$.

\section{Proof of Corollary \ref{c:cor-compact}}\label{app:cor-compact}

Let $|s|\leq R$ on $K$. For $\eta$ small, $|\eta^{1/3}s|\leq\zeta_0$ on $K$ and, by \eqref{c:l1}, $|\lambda_1(s)|\leq R+1$, so
$|e^{-\lambda_1(s)}|\leq e^{R+1}$. By \eqref{c:lemmaA},
$f_\eta(s)=e^{\lambda_1(s)}\left[\Phi(\zeta)\left(1+r_2\right)+e^{-\lambda_1(s)}r_3\right]$ with $\zeta=\eta^{1/3}s$. By
\eqref{c:Phi-bounds}, $\Phi(\zeta)=1+\eta^{2/3}s^2+O\left(\eta R^3\right)$, while
$|r_2|+|e^{-\lambda_1}r_3|\leq Ce^{R+1}e^{-\eta^{-1/3}/4}=O(\eta)$. This yields \eqref{c:fcompact1}. Equation \eqref{c:fcompact2} follows from $\lambda_1(s)+s=O(\eta)$ uniformly on $K$, by \eqref{c:compact-roots}.

\section{Proof of Corollary \ref{c:cor-resonant}}\label{app:cor-resonant}

Let $R:=\max_Q\left(|\sigma|+\psi\right)$. For $\eta\leq1$, $|s|\leq R\eta^{-1/4}$ and $|\zeta|=|\eta^{1/3}s|\leq R\eta^{1/12}$. Choose
$\eta_Q\leq\eta_A$ with $R\eta_Q^{1/12}\leq\zeta_0$. By the binomial theorem,
\[
\eta s^4=\psi^4\left(1-\frac{i\sigma\eta^{1/4}}{\psi}\right)^4
=\psi^4-4i\sigma\psi^3\eta^{1/4}-6\sigma^2\psi^2\eta^{1/2}+4i\sigma^3\psi\eta^{3/4}+\sigma^4\eta,
\]
so $|\eta s^4-\psi^4|\leq15R^4\eta^{1/4}$ and $|\eta s^4|\leq16R^4$. By \eqref{c:l1},
$\lambda_1+s+\psi^4=-\left(\eta s^4-\psi^4\right)-\eta s^4r_1$ with $|r_1|\leq C_0\eta|s|^3\leq C_0R^3\eta^{1/4}$, hence
$|\lambda_1+s+\psi^4|\leq C'\eta^{1/4}$ and $|e^{\lambda_1+s}-e^{-\psi^4}|\leq2C'\eta^{1/4}$ for $\eta$ small. By
\eqref{c:Phi-bounds}, $|\Phi(\zeta)-1|\leq|\zeta|^2+C_\Phi|\zeta|^3\leq\left(1+C_\Phi\right)R^2\eta^{1/6}$. Multiplying \eqref{c:lemmaA}
by $e^s$ and noticing $|e^s|=e^\sigma\leq e^R$,
\[
e^sf_\eta(s)-e^{-\psi^4}=\left(e^{\lambda_1+s}-e^{-\psi^4}\right)\Phi\left(1+r_2\right)
+e^{-\psi^4}\left[\left(\Phi-1\right)\left(1+r_2\right)+r_2\right]+e^sr_3,
\]
and each term is bounded by $C\eta^{1/6}$ for $\eta$ small, since $\eta^{1/4}\leq\eta^{1/6}$ and
$e^{-\eta^{-1/3}/4}\leq\eta^{1/6}$ for small $\eta$.

\section{Proof of Corollary \ref{c:cor-B1}}\label{app:cor-B1}

Apply Proposition~\ref{c:prop-B} with $\kappa_0:=\zeta_*^{-3/4}$. With $\zeta=\eta^{1/3}s$ we have
$\kappa=\eta^{-1/4}|s|^{-3/4}=|\zeta|^{-3/4}\leq\kappa_0$ and $N=\left(|s|/\eta\right)^{1/4}=\eta^{-1/3}|\zeta|^{1/4}\geq
\zeta_*^{1/4}\eta^{-1/3}$, which is $\geq N_0$ for $\eta$ small. Finally,
$\cos\theta=\Re s/|s|\geq-M/|s|\geq-M\eta^{1/3}/\zeta_*\geq-\sin\theta_1=\cos\left(\frac\pi2+\theta_1\right)$ for $\eta$ small, so
$|\theta|\leq\frac\pi2+\theta_1$. The bound $e^{-cN}$ is the claimed one.

\section{Proof of Corollary \ref{c:cor-polefree}}\label{app:cor-polefree}

Take $\eta_P(M):=\min\{\eta_A,\eta_B(M,\zeta_0)\}$ and use Proposition~\ref{c:prop-A} if $|\eta^{1/3}s|\leq\zeta_0$ and
Corollary~\ref{c:cor-B1} with $\zeta_*=\zeta_0$ if $|\eta^{1/3}s|\geq\zeta_0$.

\section{Proof of Lemma \ref{c:lem-lowerF}}\label{app:lem-lowerF}

$F\left(F-e^{-s}\right)=1$ gives $1\leq|F|\left(|F|+e^{-\Re s}\right)\leq|F|\left(|F|+e^M\right)$, so $|F|$ is at least the positive root of $x^2+e^Mx-1$.

\section{Proof of Lemma \ref{c:lem-damping}}\label{app:lem-damping}

Write $\sigma:=\Re s\geq0$, $\omega:=\Im s$ and $\zeta:=\eta^{1/3}s$. By \eqref{c:l1}, $\lambda_1=-s-\eta s^4\left(1+r_1\right)$ with
$|r_1|\leq C_0|\zeta|^3\leq C_0\zeta_0^3$, so that
\[
\Re\lambda_1=-\sigma-\Re\left(\eta s^4\left(1+r_1\right)\right).
\]
The dissipative term can increase $\Re \lambda_1$ only if $\Re (s^4)<0$, that is, if $|\arg s|\in\left(\frac\pi8,\frac{3\pi}8\right)$. Such
directions satisfy $|\omega|<\tan\left(\frac{3\pi}8\right)\,\sigma$. Fix $b>\tan\left(\frac{3\pi}8\right)$ and put $c_b:=\cos\left(4\arctan\left(1/b\right)\right)$,
which is positive because $4\arctan\left(1/b\right)<\frac\pi2$. Decreasing $\zeta_0$, we may assume $C_0\zeta_0^3\leq\frac{c_b}2\leq1$ and
$2\zeta_0^3\sqrt{1+b^2}\leq\frac12$.
\begin{itemize}[itemsep=1pt]
\item If $|\omega|\leq b\sigma$, then $|s|\leq\sqrt{1+b^2}\,\sigma$ and, since $\eta|s|^4=|\zeta|^3|s|$,
$|\eta s^4\left(1+r_1\right)|\leq2\zeta_0^3\sqrt{1+b^2}\,\sigma\leq\frac\sigma2$. Hence $\Re\lambda_1\leq-\frac\sigma2$.
\item If $|\omega|>b\sigma$, then $\big|\arg s\mp\frac\pi2\big|<\arctan\left(1/b\right)$, so $4\arg s$ lies within $4\arctan\left(1/b\right)$ of
$\pm2\pi$ and $\Re s^4\geq c_b|s|^4$. Therefore $\Re\left(\eta s^4\left(1+r_1\right)\right)\geq\eta|s|^4\left(c_b-|r_1|\right)\geq\frac{c_b}2\eta|s|^4\geq0$, and
$\Re\lambda_1\leq-\sigma$.
\end{itemize}

\section{Proof of Lemma \ref{c:lem-hgamma}}\label{app:hgamma}
 
Let $F:=e^sf_\eta(s)$. Then $h_\gamma=F^2-F$ and, by Corollary~\ref{c:cor-resonant} with $\psi=\varphi$,
$|F-e^{-\varphi^4}|\leq C_Q\eta^{1/6}$. In particular $|F|\leq2$ for $\eta$ small. Hence
$|h_\gamma-h_0|=\left|\left(F-e^{-\varphi^4}\right)\left(F+e^{-\varphi^4}-1\right)\right|\leq4C_Q\eta^{1/6}$.

\printbibliography
[heading=bibintoc,title={Bibliography}]
\end{document}